\documentclass[preprint,12pt,3p]{elsarticle}

\usepackage[utf8]{inputenc}
\usepackage[tbtags]{amsmath}
\usepackage{blindtext}
\usepackage{float}
\usepackage{amssymb}
\usepackage{amsthm}
\usepackage{nomencl}
\usepackage{commath}
\usepackage{url}
\usepackage[english]{babel}
\usepackage{csquotes}
\usepackage{color}
\usepackage{rotating}
\usepackage{graphicx}
\usepackage{tikz}
 \usetikzlibrary{shapes,arrows,spy,positioning,patterns,shapes.geometric}
\usepackage{pgfplots}
\usepackage{textcomp}
\usepackage[detect-weight=true, binary-units=true]{siunitx}
\usepackage{textcase}
\usepackage[english,ruled,vlined]{algorithm2e}
\usepackage{multirow,makecell}
\pgfplotsset{width=0.9\textwidth,compat=1.9}
\graphicspath{ {Figures/} }
\usepackage{subfigure,booktabs}

\makenomenclature

\usepackage{etoolbox}
\renewcommand\nomgroup[1]{%
    \item[\bfseries
    \ifstrequal{#1}{A}{Sets}{%
    \ifstrequal{#1}{B}{Parameters}{%
    \ifstrequal{#1}{C}{Binary decision variables: Discrete resources }{%
    \ifstrequal{#1}{D}{Binary decision variables: Continuous resources}{%
    \ifstrequal{#1}{E}{Continuous decision variables: Discrete resources}{%
    \ifstrequal{#1}{F}{Continuous decision variables: Continuous resources}{%
    \ifstrequal{#1}{G}{Continuous decision variables: Others}{
    \ifstrequal{#1}{H}{Binary decision variables: Lines}{}}}}}}}}%
]\vspace{0.1in}}
\addto\captionsenglish{\renewcommand*{\tablename}{TABLE}}

\usepackage{mathtools}

\newcommand{\etal}{\textit{et al}. }

\allowdisplaybreaks[2]

\graphicspath{ {Figures/} }
\newcommand{\bs}[1]{\boldsymbol{#1}}

  \usepackage[textsize=scriptsize]{todonotes}

\tikzstyle{startstop} = [rectangle, rounded corners, minimum width=3cm, minimum height=1cm,text centered, draw=black, fill=white]
\tikzstyle{io} = [trapezium, trapezium left angle=70, trapezium right angle=110, minimum width=3cm, minimum height=1cm, text centered, draw=black, fill=white]
\tikzstyle{process} = [rectangle, minimum width=3cm, minimum height=1cm, text centered, draw=black, fill=white]
\tikzstyle{decision} = [diamond, minimum width=3cm, minimum height=1cm, text centered, draw=black, fill=white]
\tikzstyle{arrow} = [thin,->,>=stealth]
\journal{Applied Energy}

\begin{document}
\sloppy

\begin{frontmatter}


\title{Algorithms for Optimal Topology Design, Placement, Sizing and Operation of Distributed Energy Resources in \\ Resilient Off-grid Microgrids}

\author[label1]{Sreenath Chalil Madathil\corref{cor1}}
\address[label1]{Watson Institute of Systems Excellence, Binghamton University, Binghamton, NY}
\address[label2]{Department of Industrial Engineering, Clemson University, Clemson, SC}
\address[label3]{Theoretical Division (T-5), Los Alamos National Laboratory, Los Alamos, NM}

\cortext[cor1]{Corresponding author}

\ead{schalil@binghamton.edu}

\author[label3]{Harsha Nagarajan}

\author[label3]{Russell Bent}

\author[label2]{Scott J. Mason}

\author[label2]{Sandra D. Ek\c{s}io\u{g}lu}

\author[label2]{Mowen Lu}

\begin{abstract}
Power distribution in remote communities often depends on off-grid microgrids. In order to address the reliability challenges for these microgrids, we develop a mathematical model for topology design, capacity planning, and operation of distributed energy resources in microgrids that includes N-1 security analysis. Due to the prohibitive size of the optimization problem, we develop a rolling-horizon algorithm that is combined with scenario-based decomposition to efficiently solve the model. We demonstrate the efficiency of our algorithm on an adapted IEEE test network and a real network from an Alaskan microgrid. We also compare our model's required solution time with commercial solvers and recently developed decomposition algorithms to solve similar problems.
\end{abstract}

\begin{keyword}
Off-grid microgrid \sep  Rolling-Horizon \sep N-1 security-constrained power flow \sep Decomposition algorithm \sep Topology design
\end{keyword}

\end{frontmatter}


\section{Introduction}
%
%
%
%


Within the United States and many other areas of the world, remote communities are disconnected from bulk transmission systems. Given the economic hurdles associated with connecting remote communities to these systems, many will remain isolated for the foreseeable future. However, it is important that these communities have the same level of reliability afforded by the bulk transmission systems \cite{DoeromdstRC2015}. To address this need, we develop an expansion planning model for off-grid microgrids that balances the costs of designing the system with the cost of operating these grids under N-1 reliability criteria. This model includes all three major decisions associated with the design and operation of off-grid microgrids: identifying the installation locations of power sources, determining capacity and power dispatch of those resources, and prescribing the network topology \cite{lasseter2004microgrid}. Though critically important, this problem is very difficult to solve given the non-convexities in discrete installation choices and power flow physics.

To address this problem, we use the mixed-integer, quadratically constrained, quadratic programming (MIQCQP) resource planning model of \cite{chalil2017remote,Mashayekh2016} and generalize it to support expansion planning with N-1 security constraints on lines. The resulting model is significantly more challenging to solve (the methods of \cite{chalil2017remote,Mashayekh2016} do not directly scale to this problem) and we develop a novel rolling-horizon (RH) algorithm to solve this problem.
In short, the key contributions  of this paper are:
\begin{itemize}
    \item To the best of our knowledge, this is the first  distribution  systems planning model with topology decisions and N-1 reliability on lines  that incorporates nonlinear ac physics,
    time-extended operations, distributed energy resource planning decisions (including conventional generators and batteries), and power-device efficiencies. We refer to this problem as the ac integrated resource planning problem for microgrids (ACIRPM).
    \item An algorithm which decomposes the problem by time and scenario that efficiently solves this problem.
    \item A demonstration on real system data and empirical validation of the results when compared with state-of-the-art MIQCQP solvers.
\end{itemize}

\begin{figure*}
\vspace{0.2cm}
\centering
\includegraphics[width=\textwidth]{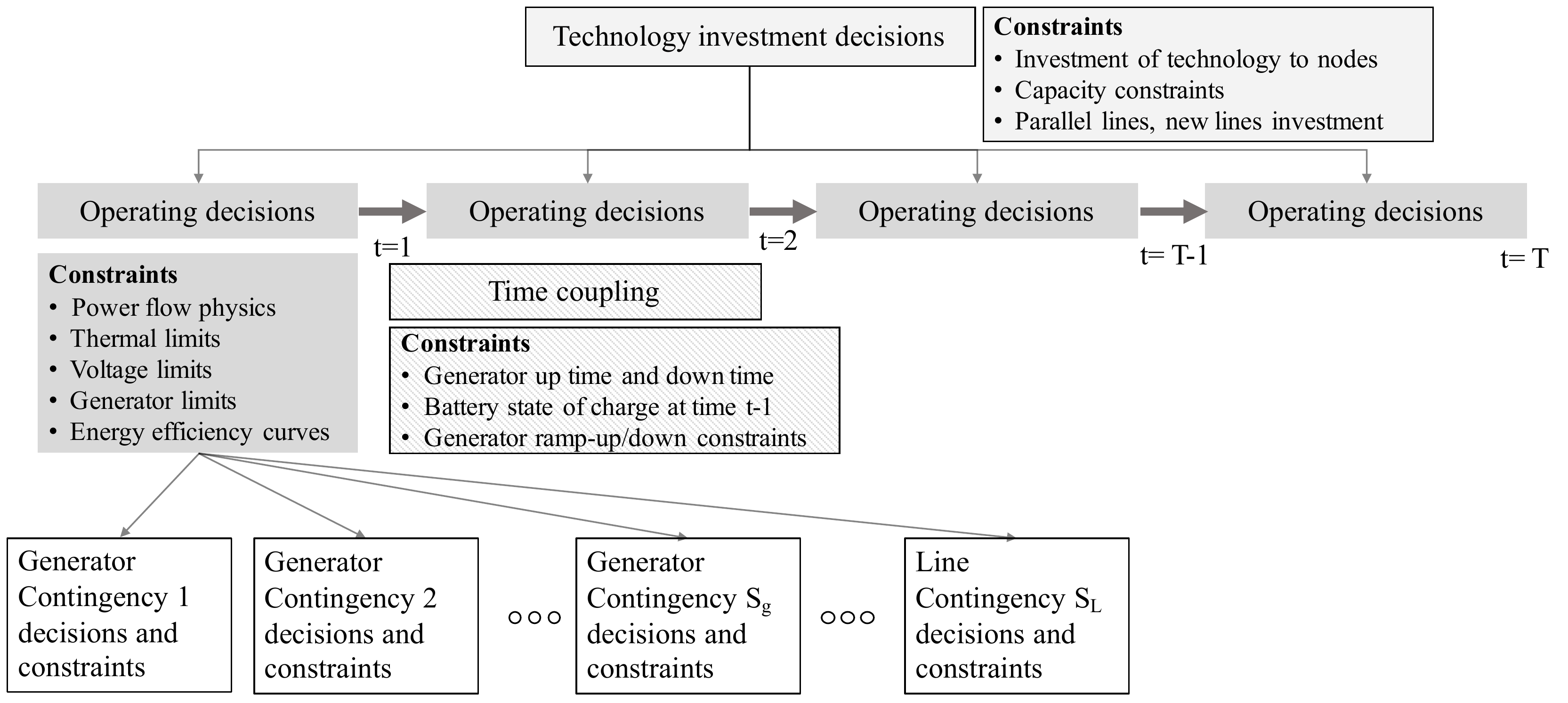}
\caption{Modified flowchart of the microgrid resource planning model described in \cite{chalil2017remote}. 
This flow chart expresses the stages and time-scales of the decisions modeled in our expansion planning model. The top level includes all planning decisions. This level includes our contribution, the inclusion of topology investment decisions. These investment decisions are used to determine time-extended operating decisions (second level). The operating decisions are constrained by N-1 security constraints (third level) and are applied at all operating points. Here, our contribution is the modeling of line contingencies.}
\label{Fig:SchematicDiag-P2}
\end{figure*}

\noindent

\paragraph{Literature Review} 
The importance of topology designs in power systems are discussed extensively in many research papers. An early study of the investment and operation of multiple energy systems along with their topology is presented in \cite{bakken2007etransport}. This paper examines the design and operation of multiple energy carriers within a locality and considers suggestions for alternate locations that satisfy pre-defined future demands. The authors recommend ranking the installation of various energy carriers over a defined time line. Their solution selects the network topology that minimizes the total cost of installation and operation. Even though \cite{bakken2007etransport} considers topological design decisions, no contingency analysis was performed on the system.  Furthermore, \cite{crucitti2004topological} studies an estimation for the vulnerability of an electric grid using topological analysis.   The authors model cascading failures of power systems via dynamic load redistribution on the network and observe that the system is highly vulnerable when heavily loaded nodes are removed from the system.

In many places in the literature, the optimization of microgrid topology design and operations is discussed as one of the main research needs in power systems, i.e. \cite{perez2016review, mashayekh2017mixed}. This observation has driven a number of studies on how topology impacts system security 
\cite{albert2004structural,crucitti2004topological,bakken2007etransport}. As noted in \cite{lasseter2002microgrid}, typical microgrid architectures 
are organized in groups of radial feeders that are part of either a distribution network or in independent, remote, off transmission locations. Under these architectures, the removal of sets of nodes (network disruptions, generator breakdowns or line failures), can lead to cascading failures.  Given this observation, \cite{albert2004structural} considers the robustness of power systems from a topological perspective and verifies correlations between reliability and redundancy of network structure and emphasizes that a redundant network enhances reliability.   In all of these papers, focus is placed on analyzing existing topology choices and these papers do not discuss how to design the topology of the network.

Some of the techniques for enhancing microgrid reliability using topological designs include interconnected microgrids and network redundancy \cite{erol2011reliable,kahveci2016optimization,zinchenko2016optimal}. Works by \cite{erol2011reliable} and \cite{kahveci2016optimization} discuss the role of interconnected microgrids.  Kahveci \etal \cite{kahveci2016optimization} present topology layouts derived from heuristics that employ clustering and graph theoretic methods. The authors discuss a heuristic approach to topology design for both ``greenfield'' sites and the expansion of existing military distribution networks. The algorithm first identifies a minimum spanning tree between various nodes and next identifies clusters that are electro-mechanically stable during islanding conditions. Unfortunately, this paper considers only topology design and did not consider the operational aspects of the microgrid, technology citing, and capacity design. Zinchenko \etal \cite{zinchenko2016optimal} solve the transmission expansion planning problem with line redundancies as a two-path problem using a variant of Dijkstra's algorithm to find the shortest path between origin and destination. They conclude that redundancy may be the only option to ensure resiliency in power systems. Hence in our paper, we consider installation of parallel, redundant lines to ensure that the network is N-1 secure for line contingencies. To the best of our knowledge, no other previous research efforts consider N-1 security analysis on line contingencies in microgrids.

The most closely related work to this paper is found in \cite{chalil2017remote,Mashayekh2016}. These papers develop a resource planning model for optimizing the placement of generation resources that enforces N-1 generator reliability on microgrids.  Though these papers also consider time-extended operations, the algorithms presented do not scale and can take up to 24 hours to solve the problem (on easier instances). Further, they do not consider expansion planning or N-1 reliability on lines.  These two modeling details, in combination with time-extended operation, can significantly increase the complexity of the problem and necessitate the development of new algorithmic approaches.
More generally, the power engineering community has developed a number of techniques for solving problems with multiple time periods like ACIRPM. These methods include Benders decomposition \cite{alguacil2000multiperiod,fortenbacher2018optimal}, rolling-horizon (RH) methods \cite{palma2013microgrid}, graph partitioning \cite{che2016optimal}, and branch-and-bound algorithms coupled with Langrangian-dual relaxation \cite{gopalakrishnan2013global}. Considering the strength of recent RH methods in industrial domains, such as supply chain optimization \cite{zamarripa2016rolling} and  microgrids \cite{palma2013microgrid}, we developed a modified RH approach for solving the ACIRPM. Uniquely, we consider different approaches, such as scenario-based decomposition (SBD), for solving the sub problems constructed by the proposed RH approaches.

The rest of the paper is organized as follows: Section II proposes a mathematical formulation for the resilient design and operation of off-grid microgrids  with N-1 security constraints on generators and lines. Section III presents a rolling-horizon algorithm for solving the model efficiently and compares its results with scenario-based decomposition method. Next, numerical results on two case studies are discussed in Section IV. Finally, Section V presents  conclusions and future research directions.

\section{Mathematical Formulation}
In this section, we present the ACIRPM model. The ACIRPM model combines expansion planning decisions with time extended operations, resource planning, efficiencies, and N-1 security criteria to optimize a microgrid for resilience.

\noindent \textit{Model parameters and variables:} In this paper, all constant parameters are typeset in \textbf{bold}.

\nomenclature[A,01]{$\mathcal{N}$}{set of nodes (buses), indexed by $i$}
\nomenclature[A,02]{$\mathcal{E}$}{set of existing edges (lines and transformers), indexed by $e$. Each edge is assigned an arbitrary direction from a bus $i$ to a bus $j$}
\nomenclature[A,03]{$\mathcal{E}_n$}{set of new edges (lines and transformers), indexed by $e$. Each edge is assigned an arbitrary direction from a bus $i$ to a bus $j$}
\nomenclature[A,07]{$\mathcal{E}_i^+$}{set of existing and new edges connected to bus $i$ and oriented from $i$, indexed by $e$}
\nomenclature[A,08]{$\mathcal{E}_i^-$}{set of existing and new edges connected to bus $i$ and oriented to $i$, indexed by $e$}
\nomenclature[A,09]{$\mathcal{T}$}{set of time periods, indexed by $t$, numbered from 1 to $\abs{\mathcal{T}}$}
\nomenclature[A,11]{$\mathcal{C}$}{set of continuous resources, indexed by $c$}
\nomenclature[A,12]{$\mathcal{C}^D\subseteq \mathcal{C}$}{set of continuous resources with discrete operation, indexed by $c$}
\nomenclature[A,13]{$\mathcal{C}^C\subseteq \mathcal{C}$}{set of continuous resources with continuous operation, indexed by $c$}
\nomenclature[A,14]{$\mathcal{C}^B\subseteq \mathcal{C}^C$}{set of continuous battery resources, indexed by $c$}
\nomenclature[A,15]{$\mathcal{D}$}{set of discrete resources, indexed by $d$}
\nomenclature[A,16]{$\mathcal{D}^D\subseteq \mathcal{D}$}{set of discrete resources with discrete operation, index by $d$}
\nomenclature[A,17]{$\mathcal{D}^C\subseteq \mathcal{D}$}{set of discrete resources with continuous operation, indexed by $d$}
\nomenclature[A,18]{$\mathcal{A} = \mathcal{C} \cup  \mathcal{D} $}{set of all resources, indexed by $a$}
\nomenclature[A,19]{${\mathcal{C}}_i^C\subseteq \mathcal{C}$}{set of continuous resources at bus $i$, indexed by $c$}
\nomenclature[A,20]{$\mathcal{C}^{CB}_i \subseteq \mathcal{C}$}{set of continuous resources with storage capabilities at bus $i$, indexed by $c$}
\nomenclature[A,21]{$\mathcal{D}_i\subseteq \mathcal{D}$}{set of discrete resources at bus $i$, indexed by $d$}
\nomenclature[A,22]{$\mathcal{A}_i\subseteq \mathcal{A}$}{set of resources at bus $i$, indexed by $a$}
\nomenclature[A,23]{$\Omega$}{set of scenarios for N-1 security analysis, indexed by $\omega$}
\nomenclature[A,24]{$\mathcal{A}_i^\omega \subseteq \mathcal{A}$}{set of resources available at bus $i$, during the contingency $\omega \in \Omega$,   indexed by $a$}
\nomenclature[A,25]{$\mathcal{A}^\omega \subseteq \mathcal{A}$}{set of resources available, during the contingency $\omega \in \Omega$,   indexed by $a$}
\nomenclature[A,26]{$\mathcal{E}^\omega \subseteq \mathcal{E}$}{set of available edges, during the contingency $\omega \in \Omega$, indexed by $e$}

\nomenclature[B,01]{$\bs{f}_a$}{fixed cost for resource $a \in \mathcal{A}$, (\$)}
\nomenclature[B,02]{$\bs{g}_a$}{variable cost for resource $a \in \mathcal{A}$, ($\$/\textrm{MW}$)}
\nomenclature[B,03]{${\bs{\kappa}_{a,0}}, {\bs{\kappa}_{a,1}}, {\bs{\kappa}_{a,2}}$ }{fixed, linear, and quadratic operational cost for resource $a \in \mathcal{A}$, ($\$$)}
\nomenclature[B,04]{${\bs{f}_{e}}$}{installation cost for line $e \in \mathcal{E}_n$, (\$)}
\nomenclature[B,05]{${\overline{\bs{u}}_d}$ , ${\underline{\bs{u}}_d}$}{minimum up-time and down-time for resource $d \in \mathcal{D}^D$, (time-step)  }
\nomenclature[B,06]{$\overline{\bs{\gamma}}_d$ , $\underline{\bs{\gamma}}_d$}{ramp up and ramp down rate for resource $d \in \mathcal{D}^D$, ($\textrm{MW}$/time-step)}
\nomenclature[B,08]{$\bs{s}_{e}$}{apparent power thermal limit on line $e \in \mathcal{E}$, ($\textrm{MVA}$)}
\nomenclature[B,10]{$\bs{dp}_i^t, \bs{dq}_i^t$}{active and reactive power demand at bus $i \in \mathcal{N}$ at time $t \in \mathcal{T}$, ($\textrm{MW}, \textrm{MVAr}$)}
\nomenclature[B,11]{$\overline{{\bs{pg}}}_{a}, \overline{{\bs{qg}}}_{a}$}{maximum active and reactive power generated by a resource $a \in \mathcal{A}$, ($\textrm{MW}, \textrm{MVAr}$)}
\nomenclature[B,12]{$\underline{{\bs{pg}}}_{a}, \underline{{\bs{qg}}}_{a}$}{minimum active and reactive power generated by a resource $a \in \mathcal{A}$, ($\textrm{MW}, \textrm{MVAr}$)}
\nomenclature[B,13]{$\bs{\delta}_{a}$}{maximum droop coefficient for  generators $a \in \mathcal{A}$}
\nomenclature[B,14]{${{\bs{\Gamma}}}_{c}$}{maximum energy storage capacity of the battery $c \in \mathcal{C}^B$, ($\textrm{MW}\textrm{-}hr$)}
\nomenclature[B,15]{$\underline{\bs{v}_i},\overline{\bs{v}_i}$}{squared voltage lower and upper bound at bus $i \in \mathcal{N}$, ($(kV)^2$)}
\nomenclature[B,16]{$\overline{\bs{s}}_a$}{maximum apparent power generated by resource $a \in \mathcal{A}$, ($\textrm{MVA}$)}
\nomenclature[B,17]{$\bs{l}_a^p$}{stand-by loss (y intercept) of a resource $a \in A$ for each piecewise function $p \in \{1,..,P\}$, ($\textrm{MW}$) }
\nomenclature[B,19]{$[\bs{\eta}^1_a \ldots \bs{\eta}^p_a]$}{vector of piecewise marginal efficiencies of maximum rated power, ($\%$)}
\nomenclature[B,22]{$\bs{R}_{e}, \bs{X}_{e}$}{resistance and reactance of line $e \in \mathcal{E}$, ($k\varOmega$)}
\nomenclature[B,23]{$\bs{\Delta t}$}{duration of a time step, ($hr$)   }
\nomenclature[B,24]{$\bs{h}_{\textrm{i}}$}{maximum number of continuous resources at bus $i$}

\nomenclature[B,25]{$\bs{k}_i$}{maximum number of discrete resources at bus $i$, indexed by $k_i$}

\nomenclature[C,01]{$x_{d}^t$ }{active/inactive status for resource $d \in \mathcal{D}^D$ at time $t \in \mathcal{T}$} 
\nomenclature[C,02]{$y_{d}^t$}{start-up status for resource $d \in \mathcal{D}^D$ at time $t \in \mathcal{T}$}
\nomenclature[C,03]{$w_{d}^t$}{shut-down status for resource $d \in \mathcal{D}^D$ at time $t \in \mathcal{T}$}
\nomenclature[C,04]{${b}_{d}$}{status indicator if discrete resource $d \in \mathcal{D}$ is built}
\nomenclature[D,01]{${b}_{c}$}{status indicator if continuous resource $c \in \mathcal{C}$ is built}
\nomenclature[E,01]{${{{pg}}}^t_{d}$ }{ac active power generation during time $t \in \mathcal{T}$ for discrete resource $d \in \mathcal{D}$, ($\textrm{MW}$)}
\nomenclature[E,02]{${{{qg}}}^t_{d}$ }{ac reactive power generation during time $t \in \mathcal{T}$ for 
discrete resource $d \in \mathcal{D}$, ($\textrm{MVAr}$)}
\nomenclature[E,03]{${\widehat{pg}}^t_{d}$ }{ac active power generation before losses during time $t \in \mathcal{T}$ for  discrete resource $d \in \mathcal{D}$, ($\textrm{MW}$)}
\nomenclature[F,01]{$\widetilde{pg}_{c}, \widetilde{qg}_{c}$}{installed maximum active and reactive power generated by a resource $c \in \mathcal{C}$, ($\textrm{MW}, \textrm{MVAr}$)}
\nomenclature[F,01.5]{$\widetilde{s}_c$}{installed maximum apparent power generated by resource $c \in \mathcal{C}$, ($\textrm{MVA}$)}
\nomenclature[F,02]{${{{pg}}}^t_{c}, {{\textrm{qg}}}^t_{c}$ }{ac apparent power generation during time $t \in \mathcal{T}$ for continuous resource $c \in \mathcal{C}$, ($\textrm{MW}, \textrm{MVAr}$)}
\nomenclature[F,03]{$\widehat{pg}^t_{c}$ }{ac active power generation before losses during time $t \in \mathcal{T}$ for continuous resource $c \in \mathcal{C}$, ($\textrm{MW}$)}
\nomenclature[F,10]{$\text{\c{e}}_{c}^t$}{energy stored (state-of-charge) in the continuous resource battery $c \in \mathcal{C}^B$ at time $t \in \mathcal{T}$ , ($\textrm{MW}\textrm{-}hr$)}
\nomenclature[G,02]{${{p}}_{e}^t , {{q}}_{e}^t$}{active and reactive power flow though edge $e \in \mathcal{E}$ at time $t \in \mathcal{T}$, ($\textrm{MW}, \textrm{MVAr}$)}
\nomenclature[G,03]{$v_i^t$}{squared voltage at node $i \in \mathcal{N}$ at time $t \in \mathcal{T}$, ($(kV)^2$)}
\nomenclature[G,04]{$\textit{lp}_{i}^{t,\omega}, \textit{lq}_{i}^{t,\omega}$}{apparent power slack at node $i \in \mathcal{N}$ at time $t \in \mathcal{T}$ during contingency scenario $\omega \in \Omega$, ($\textrm{MW}, \textrm{MVAr}$)}
\nomenclature[H,01]{$b_{e}$}{status indicator if line $e\in \mathcal{E}_n$ is built}
{\small \printnomenclature[1in]}

\paragraph{Objective function}

The objective function of the ACIRPM lexicographically minimizes load slack during the contingencies and then minimizes the total installation and operation cost of energy resources and the cost of installing new lines to enhance network resiliency \eqref{P2-Objfun01}. 
\begin{subequations}
\begin{align}
\begin{split}
&\min \bigg\langle \;\;
\bigg(\sum_{i \in \mathcal{N}}\sum_{t \in \mathcal{T}}\sum_{\omega \in \Omega}(|\textit{lp}_{i,t}^\omega|+|\textit{lq}_{i,t}^\omega|)\bigg), \\
&\bigg(\sum_{c \in \mathcal{C}}\bs{f}_c{b}_{c} + \sum_{c \in \mathcal{C} \setminus \mathcal{C}^B}\bs{g}_c{\widetilde{pg}}_{c} +\sum_{c \in  \mathcal{C}^B}\bs{g}_c{\widetilde{s}}_{c}\bigg) +
\sum_{d \in \mathcal{D}}\bs{f}_d b_{d} \; +  \sum_{e \in \mathcal{E}_n}\bs{f}_{e}b_{e} +\\
&\bigg(\sum_{t \in T} \sum_{a \in \mathcal{A}}\big((\bs{\kappa}_{a,2})({\widehat{pg}}_{a}^t)^2   +(\bs{\kappa}_{a,1})({\widehat{pg}}_{a}^t)   + (\bs{\kappa}_{a,0})(b_{a})\big)\bigg) 
\;\;\bigg\rangle \label{P2-Objfun01}
\end{split}
\end{align}
\end{subequations}

\paragraph{Resource Planning}
The constraints associated with the availability of resources are defined in equations \eqref{P2-ct-01}-\eqref{P2-dt-04}. Here, equations \eqref{P2-ct-01}-\eqref{P2-ct-02} link the installed capacity of continuous resources to the build variables (${b}_{c}$).
Equation \eqref{P2-ct-02a} links the installed apparent power capacity of storage devices with the build variable.  
Similarly, equations \eqref{P2-dt-03}-\eqref{P2-dt-04} constrain the capacity limits for discrete resources. The number of continuous and discrete resources installed at a bus is constrained by equations \eqref{P2-ct-02a1}-\eqref{P2-ct-02a2}. 


\begin{subequations}
\begin{align}
&0 \leq {pg}_{c}^t \le {\widetilde{pg}}_{c} \le {b}_{c}{\bs{\overline{pg}}_c} &&\forall~ c \in \mathcal{C} \setminus \mathcal{C}^B, t \in \mathcal{T} \label{P2-ct-01}\\
&{b}_{c}\bs{\underline{qg}}_c \le {qg}_{c}^t \le {\widetilde{qg}}_{c} \le {b}_{c}\bs{\overline{qg}}_c&&\forall~ c \in \mathcal{C} \setminus \mathcal{C}^B, t \in \mathcal{T} \label{P2-ct-02}\\
&{\widetilde{s}}_{c} \le b_{c}\overline{\bs{s}}_c &&\forall~c \in \mathcal{C}^B \label{P2-ct-02a}\\
&\sum_{c \in \mathcal{C}_i}{b}_{c} \le \bs{h}_{i} &&\forall~ i \in \mathcal{N}\label{P2-ct-02a1}\\
&\sum_{d \in \mathcal{D}_i}{b}_{d} \le \bs{k}_{i} &&\forall~ i \in \mathcal{N}\label{P2-ct-02a2}\\
&0 \le {{\widehat{pg}}}_{d}^t \le \bs{\overline{{pg}}}_{d} {b}_{d} &&\forall~ d \in \mathcal{D}^C, t \in \mathcal{T} \label{P2-dt-03}\\
&\underline{\bs{qg}}_{d} {b}_{d} \le {{\widehat{qg}}}_{d}^t \le \bs{\overline{{qg}}}_{d} b_{d} &&\forall~ d \in \mathcal{D}^C, t \in \mathcal{T} \label{P2-dt-04}
\end{align}
\end{subequations}

\paragraph{Power Flow Physics}
The physics of the ACIRPM are shown in equations \eqref{P2-powerflow-01}-\eqref{P2-lc-02}, where the \textit{LinDistFlow} equations \eqref{P2-lc-02} of \cite{gan2015exact,baran1989network} are used. Here, equations \eqref{P2-powerflow-01}-\eqref{P2-powerflow-02} model Kirchoff's Law and equations \eqref{P2-lc-02} model Ohm's Law.

\begin{subequations}
\begin{align}
&\sum_{a \in \mathcal{A}_i}{pg}_{a}^t  - \bs{{dp}}_{i}^t =  \sum_{ e \in \mathcal{E}_i^+}{p}_{e}^t - \sum_{ e \in \mathcal{E}_i^-}{p}_{e}^t &&\forall~ i \in \mathcal{N}, t \in \mathcal{T} \label{P2-powerflow-01}\\
&\sum_{a \in \mathcal{A}_i}qg_{a}^t  - \bs{dq}_{i}^t =  \sum_{ e \in \mathcal{E}_i^+}{q}_{e}^t - \sum_{ e \in \mathcal{E}_i^-}{q}_{e}^t &&\forall~ i \in \mathcal{N}, t \in \mathcal{T} \label{P2-powerflow-02}\\
&v_j^t = v_i^t - 2(\bs{\bs{R}}_{e}{p}_{e}^t+\bs{\bs{X}}_{e}{q}_{e}^t)  &&\forall~ e \in \mathcal{E}, t \in \mathcal{T} \label{P2-lc-02}
\end{align}
\end{subequations}

\paragraph{Physical Limits}
The physical limits of the ACIRPM are shown in equations \eqref{P2-lc-01}-\eqref{P2-lc-03}. Equation \eqref{P2-lc-01} places thermal limits on lines and equation \eqref{P2-lc-03} places voltage magnitude limits on buses.

\begin{subequations}
\begin{align}
&({{p}}_{e}^t)^2 + ({{q}}_{e}^t)^2 \le (\bs{s}_{e})^2  &&\forall~ e \in \mathcal{E}, t \in \mathcal{T} \label{P2-lc-01}\\
&\underline{\bs{v}_i} \le {v_i^t} \le \overline{\bs{v}_i}  &&\forall~ i \in \mathcal{N}, t \in \mathcal{T} \label{P2-lc-03}
\end{align}
\end{subequations}

\paragraph{Generator Limits}
Equations \eqref{P2-dtdo-01}-\eqref{P2-dtdo-09} model the operating limits of resources defined as discrete  generators (i.e., diesel generators). The connection between a generator's on/off status and its start up and shut down time are modeled with equations \eqref{P2-dtdo-01}-\eqref{P2-dtdo-03} . Equations \eqref{P2-dtdo-04}-\eqref{P2-dtdo-05} link active and reactive power dispatch with the generator's status. Generator operating characteristics like minimum up-time, minimum down-time, ramp-up time, and ramp-down time are constrained using equations \eqref{P2-li-07}-\eqref{P2-li-08}. We model the boundary conditions for uptime and downtime using

\[
\alpha_d = \{\rho \in \mathcal{T} : t-\overline{\bs{u}}_d + 1 \le \rho \le t\} 
\]

\noindent and

\[
\zeta_d =  \{\rho \in \mathcal{T} : t-\underline{\bs{u}}_d + 1 \le \rho \le t\} 
\]

\noindent respectively.


\begin{subequations}
\begin{align}
&x_{d}^t \le b_{d}   &&\forall~ d \in \mathcal{D}^D, t \in \mathcal{T} \label{P2-dtdo-01}\\
&x_{d}^t = x_{d}^{t-1}+y_{d}^t -w_{d}^t &&\forall~ d \in \mathcal{D}^D, t \in \mathcal{T}\label{P2-dtdo-02}\\
&y_{d}^t + w_{d}^t \le 1  &&\forall~ d \in \mathcal{D}^D, t \in \mathcal{T}\label{P2-dtdo-03}\\
&0 \le {{\widehat{pg}}}_{d}^t \le \overline{\bs{pg}}_{d} x_{d}^t  &&\forall~ d \in \mathcal{D}^D,t \in \mathcal{T}\label{P2-dtdo-04}\\
&\underline{\bs{qg}}_{d} x_{d}^t \le {{\widehat{qg}}}_{d}^t \le \overline{\bs{qg}}_{d} x_{d}^t  &&\forall~ d \in \mathcal{D}^D,t \in \mathcal{T}\label{P2-dtdo-05}\\
&\sum_{\rho \in \alpha_d} y_{d}^\rho \leq x_{d}^t &&\forall~ d \in \mathcal{D}^D, t \in \mathcal{T} \label{P2-li-07} \\
&\sum_{\rho \in \zeta_d}w_{d}^\rho \leq 1-x_{d}^t &&\forall~ d \in \mathcal{D}^D, t \in \mathcal{T} \label{P2-li-08} \\
&\overline{\bs{\gamma}}_d \ge {pg}_{d}^{t} - {pg}_{d}^{t-1} -\overline{\bs{pg}}_{d} y_{d}^t &&\forall~ d \in \mathcal{D}^D t \in \mathcal{T}\label{P2-dtdo-08}\\
&\underline{\bs{\gamma}}_d \ge {pg}_{d}^{t-1} - {pg}_{d}^{t} - \overline{\bs{pg}}_{d} w_{d}^t &&\forall~ d \in \mathcal{D}^D, t \in \mathcal{T}\label{P2-dtdo-09}
\end{align}
\end{subequations}

\paragraph{Battery Limits}
Equations \eqref{P2-bat-01}-\eqref{P2-bat-04} model the operating limits of resources defined as batteries. Equation \eqref{P2-bat-01} constrains the apparent power of batteries. The charging and discharging of batteries are modeled using equations \eqref{P2-bat-03}.

\begin{subequations}
\begin{align}
&({pg}_{c}^t)^2 + (qg_{c}^t)^2 \le ({\widetilde{s}}_{c})^2  &&\forall~  c \in \mathcal{C}^B, t \in \mathcal{T}\label{P2-bat-01}\\
&\text{\c{e}}_{c}^t = \text{\c{e}}_{c}^{t-1}-{\widehat{pg}}_{c}^{t} \bs{\Delta t} &&\forall~ c \in \mathcal{C}^B, t \in \mathcal{T}\label{P2-bat-03}\\
&0 \le \widetilde{s}_{c} \le \overline{\bs{s}}_c  &&\forall~ c \in \mathcal{C}^B \label{P2-bat-04}
\end{align}
\end{subequations}

\paragraph{Efficiencies}
Component efficiencies, typically represented with piecewise linear functions (indexed by $p$) and modeled using separate binary variables, are computationally challenging. To circumvent this issue, we apply convex relaxations by relaxing the feasible space of component efficiencies to a halfspace (for details, see Figure 2 in \cite{chalil2017remote}) as shown in \eqref{P2-eff01fun} - \eqref{P2-eff03fun}. 

\begin{subequations}
\begin{align}
&pg_{c}^t \le \bs{\eta}^p_c \widehat{pg}_{c}^t + b_{c}\bs{l}^p_c &&\forall~ c \in \mathcal{C}, t \in \mathcal{T}, p \label{P2-eff01fun}\\
&{pg}_{d}^t \le \bs{\eta}^p_d \widehat{pg}_{d}^t + x_{d}^t\bs{l}^p_{d}  &&\forall~ d \in \mathcal{D}^C, t \in \mathcal{T}, p\label{P2-eff02fun} \\
&{pg}_{d}^t \le \bs{\eta}^p_d \widehat{pg}_{d}^t + {b}_{d}\bs{l}^p_{d} &&\forall~ d \in \mathcal{D}^D, t \in \mathcal{T}, p \label{P2-eff03fun}
\end{align}
\end{subequations}

\paragraph{Expansion Planning}
On/off constraints are used to model thermal limits \eqref{P2-pf-01} and Ohm's laws \eqref{P2-pf-02}-\eqref{P2-pf-03} for new lines. 
Here $\bs{M}_i= \overline{\bs{v}_i} - \underline{\bs{v}_i}$.

\begin{subequations}
\begin{align}
&({p}_{e}^t)^2 + ({q}_{e}^t)^2 \le b_{e} \cdot (\bs{s}_{e})^2 && \forall e \in {\cal{E}}_n, t \in \mathcal{T}  \label{P2-pf-01}\\
&v_j^t - v_i^t \ge -2(\bs{R}_{e}{p}_{e}^t+\bs{X}_{e}{q}_{e}^t) - \bs{M}_i(1-b_{e}) && \forall e \in {\cal{E}}_n, t \in \mathcal{T}   \label{P2-pf-02}\\
&v_j^t - v_i^t \le -2(\bs{R}_{e}{p}_{e}^t+\bs{X}_{e}{q}_{e}^t) + \bs{M}_i(1-b_{e}) && \forall e \in {\cal{E}}_n, t \in \mathcal{T}   \label{P2-pf-03}
\end{align}
\end{subequations}

\paragraph{Generator Contingencies}
Each generator contingency replicates equations \eqref{P2-ct-01}-\eqref{P2-pf-03} on subsets of $\mathcal{C}$ and $\mathcal{D}$ \cite{chalil2017remote}. The subsets remove the generators that are outaged in the contingency. Equations \eqref{P2-powerflow-01} and \eqref{P2-powerflow-02} are replaced with their load slack equivalents. During generator contingencies, note that $\mathcal{E}^\omega = \mathcal{E}$ since all the lines are available. The constraints are summarized below:

\begin{subequations}
\begin{align}
&{lp}_{i}^{t,\omega} + \sum_{a \in \mathcal{A}_i^\omega}{pg}_{a}^{t,\omega}  - \bs{dp}_{i}^{t}  =  \sum_{ e \in \mathcal{E}_i^{\omega,+}}{p}_{e}^{t,\omega} - \sum_{ e \in \mathcal{E}_i^{\omega,-}}{p}_{e}^{t,\omega} &&\forall~ i \in \mathcal{N}, t \in \mathcal{T}, \omega \in \Omega \label{P2-powerflow-01-ls}\\
& lq_{i}^{t,\omega} + \sum_{a \in \mathcal{A}_i^\omega}qg_{a}^{t,\omega}  - \bs{dq}_{i}^{t} =  \sum_{ e \in \mathcal{E}_i^{\omega,+}}{q}_{e}^{t,\omega} - \sum_{ e \in \mathcal{E}_i^{\omega,-}}{q}_{e}^{t,\omega} &&\forall~ i \in \mathcal{N}, t \in \mathcal{T}, \omega \in \Omega \label{P2-powerflow-02-ls}\\
&({{p}}_{e}^{t,\omega})^2 + ({{q}}_{e}^{t,\omega})^2 \le (\bs{s}_{e})^2  &&\forall~ e \in \mathcal{E}^\omega, t \in \mathcal{T}, \omega \in \Omega  \label{P2-lc-01-}\\
&v_j^{t,\omega} = v_i^{t,\omega} - 2(\bs{R}_{e}{{p}}_{e}^{t,\omega}+\bs{X}_{e}{{q}}_{e}^{t,\omega})  &&\forall~ e \in \mathcal{E}^\omega, t \in \mathcal{T}, \omega \in \Omega \label{P2-lc-02-ls}\\
&\underline{\bs{v}_i} \le {v_i^{t,\omega}} \le \overline{\bs{v}_i}  &&\forall~ i \in \mathcal{N}, t \in \mathcal{T}, \omega \in \Omega \label{P2-lc-03-ls}\\
&{pg}_{a}^{t} - \bs{\delta}_a x_a^t \le {pg}_{a}^{t,\omega}  \le {pg}_{a}^{t} + \bs{\delta}_a x_a^t &&\forall~  a \in \mathcal{A}^\omega, t \in \mathcal{T}, \omega \in \Omega \label{P2-droop-01-ls}\\ 
&{qg}_{a}^{t} - \bs{\delta}_a x_a^t \le {qg}_{a}^{t,\omega}  \le {qg}_{a}^{t} + \bs{\delta}_a x_a^t &&\forall~  a \in \mathcal{A}^\omega, t \in \mathcal{T}, \omega \in \Omega
\label{P2-droop-02-ls}
\end{align}
\end{subequations}

\paragraph{Line Contingencies}
Each line contingency replicates equations \eqref{P2-powerflow-01-ls}-\eqref{P2-droop-02-ls} on subsets of $\mathcal{E}$. The subsets remove the lines that are outaged in the contingency. 

For notational simplicity, we refer to the entire model defined by \eqref{P2-Objfun01}-\eqref{P2-powerflow-02-ls} as $\mathcal{M}$. A special case of $\mathcal{M}$ relaxes the $N-1$ constraints to a subset $\hat{\Omega} \in \Omega$, and $\mathcal{M}_{\hat{\Omega}}$ is used to denote this special case. Similarly, a special case of $\mathcal{M}$ relaxes $\mathcal{M}$ to a subset of time points $\tau \in \mathcal{T}$. We use 
$\mathcal{M}^{\tau}$ to denote this special case.
We note that in some applications discrete resources are restricted to \textit{slots} at a node due to sizing requirements. Without loss of generality, $\mathcal{M}$ can be modified to include slots (see \cite{chalil2017remote}).

\section{Algorithms}
\subsection{Base algorithm}

We define the base algorithm as an approach that formalizes the entire model as a single input to a state-of-the-art MIQCQP commercial solver (Gurobi V8.0). We use this approach as a comparison point.

\subsection{Scenario-based decomposition (SBD) algorithm} The SBD algorithm was first applied to energy infrastructure resiliency problems in \cite{chalil2017remote,nagarajan2016optimal,byeon2018communication}, where it was shown to have considerable computational advantages. For these problems, SBD converges to the global optimal and we use it as another comparison point for our approach. For completeness, the SBD algorithm is outlined in Algorithm \ref{P2-Alg:SBD}. SBD first relaxes the N-1 contingency constraints (model $\mathcal{M}_\emptyset$). Each N-1 contingency (scenario) is then solved given the resource and expansion planning decisions of the solution to $\mathcal{M}_\emptyset$. The constraints of the contingency that require the most load slack are then added as part of the constraint set ($S_{obj}(\omega) = |\textit{lp}_{i,t}^\omega|+|\textit{lq}_{i,t}^\omega|$).  This constraint set is added to ensure that the solution found satisfies demand under the worst contingency scenario. The SBD algorithm terminates when the load slack of the remaining contingencies are negligible or, in worst-case, all contingencies are added to the constraint set.

 \begin{algorithm}[htp]
 \caption{Scenario-based decomposition}
 \SetAlgoLined
  Create scenario set $S$, indexed by $\omega$, with all N-1 scenarios\;
  Define $S_{obj}$ as a vector of size $|S|$\; 
  \While{$\textrm{max}(S_{obj}) > 0 \textrm{ or } S \ne \emptyset$}{
     Solve the model, $\mathcal{M}_{\Omega \setminus S}$\;
     Get the values of base model decision variables, $\overline{x}$\;
     \For{$\omega \in S$}{
     Solve sub-problem for scenario $\omega$ using $\overline{x}$\;
     Update $S_{obj}(\omega)$\;
     }
     Select scenario, $\omega = \arg\max_S S_{obj}(\omega), \omega \in S$\; 
     Update scenario set  $S = S \setminus \omega$\;
     Set $S_{obj}(\omega) =  0$\;
  }
 \label{P2-Alg:SBD}
 \end{algorithm}

\subsection{Rolling-horizon (RH) algorithm}
In our initial computational experiments we found that the ACIRPM  was computationally very challenging to solve with exact methods. Here, we discuss our RH heuristic which decomposes the ACIRPM  into a sequence of smaller problems. Each of these sub-problems considers a limited number of time periods. Each sub-problem of the RH is solved using one of these above exact methods (base algorithm/SBD). Note that, solving sub-problems using exact methods does not necessarily imply that the solution of RH algorithm converges to the global optimum of the given full problem.

\begin{figure}[ht!]
    \centering
    \begin{tikzpicture}[scale=1.0]
        \draw[->,gray, thick] (0,0) -- (11,0); 
        \draw[gray, thick] (0,0) -- (0,8.6); 
        \draw[gray, dashed] (0,7.3) -- (10.5,7.3); 
        \draw[gray, dashed] (0,5.3) -- (10.5,5.3); 
        \draw[gray, dashed] (0,2.9) -- (10.5,2.9); 
        \draw[gray, dashed] (0,0.3) -- (10.5,0.3); 
        \node[] at (9,-0.5) {Design horizon};
        \node[rotate=90] at (-0.5,7) {Iterations};
        \draw[pattern color=gray!40,pattern=dots] (0.2,8.5) rectangle (9.2,7.5);
        \node[] at (4.6,8) {Schedule horizon, $\mathcal{T}$};
        \draw[pattern color=orange!50,pattern= horizontal lines] (0.2,6.7) rectangle (5.2,5.7);
        \node[] at (2.6,6.2) {Prediction horizon, $\mathcal{T}^p_1$};
        \node[rotate=90] at (10,6.3) {iteration 1};
        \path[draw=black!20,solid,line width=2mm,fill=black!20,
preaction={-triangle 90,thin,draw=black!20,shorten >=-1mm}
] (3.3, 5.4) -- (3.3, 4.2);
        \node[rotate=0] at (5.5,4.6) {Warm-start from $\mathcal{T}^p_1$};
        \draw[pattern color=purple!40,pattern=crosshatch dots] (0.2,4.2) rectangle (2.3,5.2);
        \node[] at (1.4,4.7) {$\mathcal{T}^c_{1}$};
        \draw[pattern color=orange!50,pattern= horizontal lines] (2.2,4.0) rectangle (7.2,3.0);
        \node[] at (4.7,3.5) {Prediction horizon, $\mathcal{T}^p_2$};
        \node[rotate=90] at (10,4.1) {iteration 2};
        \path[draw=black!20,solid,line width=2mm,fill=black!20,
preaction={-triangle 90,thin,draw=black!20,shorten >=-1mm}
] (5.2, 2.9) -- (5.2, 1.7);
        \node[rotate=0] at (7.4,2.1) {Warm-start from $\mathcal{T}^p_2$};
        \draw[pattern color=purple!40,pattern=crosshatch dots] (2.3,1.8) rectangle (4.4,2.8);
        \node[] at (3.45,2.3) {$\mathcal{T}^c_{2}$};
        \draw[pattern color=orange!50,pattern= horizontal lines] (4.2,1.5) rectangle (9.2,0.5);
        \node[] at (6.7,1.0) {Prediction horizon, $\mathcal{T}^p_3$};
        \node[rotate=90] at (10,1.4) {iteration 3};
    \end{tikzpicture}    
    \caption{Schematic diagram for rolling-horizon}
    \label{fig:RH}
\end{figure}
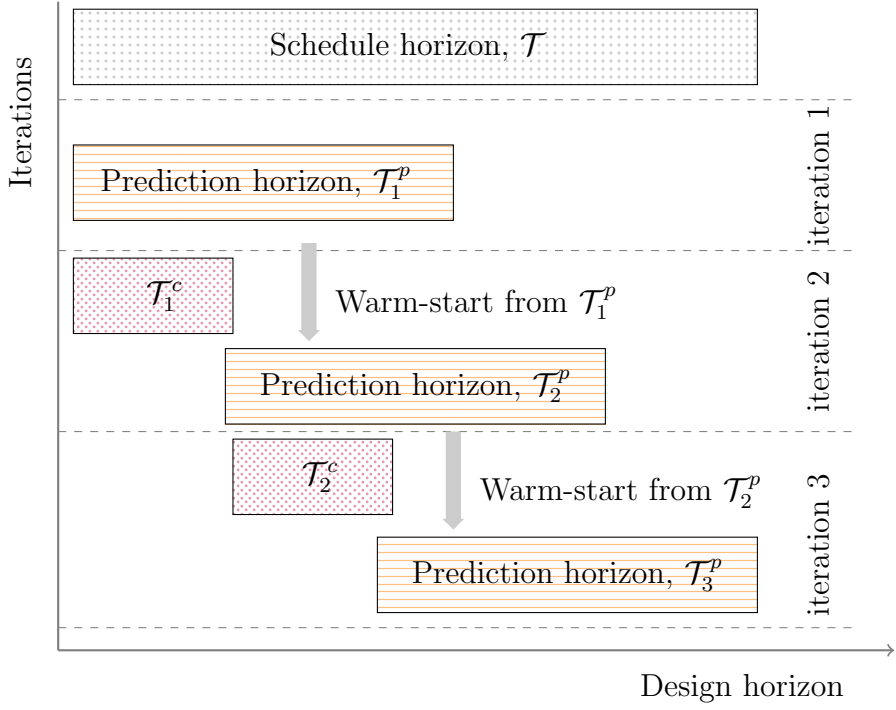

The RH algorithm is defined by three parameters, the scheduling horizon $\mathcal{T}^s=|\mathcal{T}|$, a prediction horizon $\mathcal{T}^p$, and a control horizon $\mathcal{T}^c$. The scheduling horizon defines the full length of the ACIRPM. The prediction horizon controls the size of the sub-problems that are solved, and the control horizon determines how much of the sub-problem solution is executed. More formally, let $\sigma_\tau$ denote the solution to an ACIRPM  starting at time $\tau$ and let $\sigma_\tau(\cdot)$ denote the variable assignment of $\cdot$ in solution $\sigma_\tau$. We can then recursively define the problem, $\mathcal{M}^\tau$, as $\mathcal{M}$ where $\mathcal{T}^\tau=\{t \in \mathcal{T} : t \le \tau+\mathcal{T}^p\}$ and 
extra constraints \eqref{eq:built_boundary}-\eqref{P2-li-11}.
We also define $\tilde{\tau}^c = \tau-\mathcal{T}^c$, to denote the starting time for a previous iteration's control horizon.%
\begin{subequations}\label{eq:built_boundary}
\begin{align}
 &{b}_{c} \ge \sigma_{\tilde{\tau}^c}({b}_{c})  &&\forall ~  c \in \mathcal{C} \label{P2-li-01}\\
 &{b}_{d} \ge \sigma_{\tilde{\tau}^c}({b}_{d})  &&\forall ~  d \in \mathcal{D}^D  \label{P2-li-02}\\
 &{b}_{e} \ge  \sigma_{\tilde{\tau}^c}({b}_{e}) &&\forall ~ e \in \mathcal{E} \cup \mathcal{E}_n\label{P2-li-03}\\
 &{\widetilde{{{pg}}}_{c}} \geq \sigma_{\tilde{\tau}^c}(\widetilde{{{pg}}}_{c})~~, ~~
 \widetilde{{{qg}}}_{c} \geq \sigma_{\tilde{\tau}^c}(\widetilde{{{qg}}}_{c}) &&\forall ~  c \in \mathcal{C} 
\end{align}
\end{subequations}

\noindent
Equations \eqref{eq:built_boundary} are used to enforce consistency of installation decisions between $\mathcal{T}^c$ problems.  Similarly, we also add constraints that enforce consistency in operation between $\mathcal{T}^c$ problems:

\begin{subequations}
\begin{align}
&x_{d}^\tau = \sigma_{\tilde{\tau}^c}(x_{d}^{\tau-1})+y_{d}^\tau -w_{d}^\tau &&~\forall~  d \in \mathcal{D}^D \label{P2-li-06}\\
&\sum_{\rho \in \alpha_d} y_{d}^\rho \leq x_{d}^t &&~\forall~ d \in \mathcal{D}^D, t \in \mathcal{T}^\tau \label{P2-li-07-1} \\
&\sum_{\rho \in \zeta_d}w_{d}^\rho \leq 1-x_{d}^t &&~\forall~ d \in \mathcal{D}^D, t \in \mathcal{T}^\tau \label{P2-li-08-1} \\
&x_{d}^t = \begin{cases}
    1, & \forall~ t \in [\tau,\tau+\overline{\bs{u}}_d-\widetilde{t}] ~ \text{, if } \sigma_{\tilde{\tau}^c}(y_{d}^{\widetilde{t}}) = 1\\
    0,              & \forall~ t \in [\tau,\tau+\underline{\bs{u}}_d-\widetilde{t}] ~\text{, if } \sigma_{\tilde{\tau}^c}(w_{d}^{\widetilde{t}}) = 1
\end{cases}\label{P2-yw-01}\\
&\overline{\gamma}_d \ge \sigma_{\tilde{\tau}^c}({pg}_{d}^{\tau-1}) - {pg}_{d}^{\tau} &&~\forall~  d \in \mathcal{D}^D \label{P2-li-09} \\
&\underline{\gamma}_d \ge {pg}_{d}^\tau - \sigma_{\tilde{\tau}^c}({pg}_{d}^{\tau-1}) &&~\forall~  d \in \mathcal{D}^D\label{P2-li-10}
\end{align}
\end{subequations}

Constraints \eqref{P2-li-06} use generator on/off status in $\sigma_{\tilde{\tau}^c}$ as a boundary condition. The minimum value for the up-time and downtime is updated as $\zeta_d=\max(t-\overline{\bs{u}}_d + 1, \min ({\mathcal{T}^p}))$ and $\alpha_d=\max(t-\underline{\bs{u}}_d + 1, \min ({\mathcal{T}^p}))$ respectively, for use in constraints \eqref{P2-li-07-1} and \eqref{P2-li-08-1}. The boundary condition ensures the proper calculation of minimum generator up-time and down-time across the whole planning horizon.
Finally, constraints \eqref{P2-li-09} and \eqref{P2-li-10} link the ramp-up and ramp-down rates between two adjacent time steps.


Constraints are also added that enforce consistency in operation of batteries in $\mathcal{T}^c$ problems.
Constraints \eqref{P2-li-11} ensure that the value of charge is carried forward from previous iterations, starting from time step $\mathcal{\tilde{T}}^c$.%
\begin{subequations}
\begin{align}
&\text{\c{e}}_{c}^\tau = \sigma_{\tilde{\tau}^c}(\text{\c{e}}_{c}^{\tau-1})-{pg}_{c}^\tau \bs{\Delta t} &&\forall ~ c \in \mathcal{C}^{CB}_i,   \label{P2-li-11}
\end{align}
\end{subequations}

The pseudo-code for our RH is given in Algorithm \ref{Alg:RH} and a schematic diagram of the algorithm is presented in Fig. \ref{fig:RH}. 
Each iteration uses the solution $\sigma_{\tilde{\tau}^c}$ to \textit{warm-start} the MIQCQP solver used to solve $\mathcal{M}^\tau$. The \textit{warm-start} initializes the assignment of variables in $\mathcal{M}^\tau$ with the assignments of those variables in $\sigma_{\tilde{\tau}^c}$ (where there is overlap between current and previous iterations).
 
\begin{algorithm}[ht!]
\caption{Rolling-horizon algorithm}
\SetAlgoLined
 \While{$\tau \le T$}{
    Warm-start $\mathcal{M}^\tau$ with $\sigma_{\mathcal{\tilde{T}}^c}$\;
    $\sigma_{\tau} \leftarrow$ Solve $\mathcal{M}^\tau$\;
    $\tau \leftarrow \tau + \mathcal{T}^c$\;
 }
\label{Alg:RH}
\end{algorithm}

\section{Numerical Results}
The numerical results were performed using a Microsoft Windows\textsuperscript{\textregistered} 10 Enterprise 2017 with an Intel\textsuperscript{\textregistered} Xenon\textsuperscript{TM} E5-2620 CPU @ $\SI{2.00}{\giga\hertz}$  processor with 6 cores and $\SI{56}{\giga\byte}$ RAM. The algorithms are modeled using JuMP in Julia \cite{DunningHuchetteLubin2017} and use Gurobi V8.0 \cite{gurobi} to solve the MIQCQPs. We test the performance of the algorithm and validate the model on an adapted version of the IEEE 13 node test feeder \cite{kersting2001radial} and a real microgrid from Alaska \cite{nome_data}.


\subsection{Case study 1: IEEE 13 node test feeder}

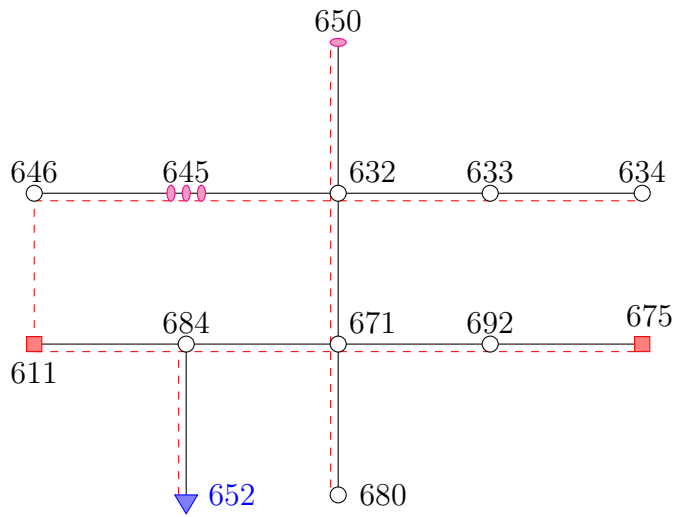
\begin{figure}[ht!]
\centering
\begin{tabular}{cc}
\begin{tikzpicture}[scale=1]
\draw (0,2) -- (8,2);
\draw (0,4) -- (8,4);
\draw (2,0) -- (2,2);
\draw (4,0) -- (4,6);
\draw[red,dashed] (0,1.9) -- (8,1.9);
\draw[red,dashed] (0,3.9) -- (8,3.9);
\draw[red,dashed] (1.9,0) -- (1.9,1.9);
\draw[red,dashed] (3.9,0.1) -- (3.9,6);
\draw[red,dashed] (0,2) -- (0,3.9);
\draw [fill=white] (0,4) circle (3.0pt) node[anchor=south] {646} ;
\filldraw[fill=magenta!50,draw=magenta] (2,4) ellipse (1.5pt and 3pt) node[anchor=south] {645};
\filldraw[fill=magenta!50,draw=magenta] (1.8,4) ellipse (1.5pt and 3pt);
\filldraw[fill=magenta!50,draw=magenta] (2.2,4) ellipse (1.5pt and 3pt);
\draw [fill=white] (4,4) circle (3.0pt) node[anchor=south west] {632};
\draw [fill=white] (6,4) circle (3.0pt) node[anchor=south] {633};
\draw [fill=white] (8,4) circle (3.0pt) node[anchor=south] {634};
\filldraw[fill=magenta!50,draw=magenta] (4,6) ellipse (3pt and 1.5pt) node[anchor=south] {650};
\filldraw[fill=red!50, draw=red] (-0.1,2.1) rectangle (0.1,1.9) node[anchor=north] at (0,1.9) {611};
\draw [fill=white] (2,2) circle (3.0pt) node[anchor=south] {684};
\draw [fill=white] (4,2) circle (3.0pt) node[anchor=south west] {671}; 
\draw [fill=white] (6,2) circle (3.0pt) node[anchor=south] {692};
\filldraw[fill=red!50, draw=red]  (7.9,1.9) rectangle (8.1,2.1)  node[anchor=south] {675};
\draw [fill=blue!50,draw=blue] (1.85,0)--(2.15,0)--(2,-0.25)--cycle;
\node [blue, fill=white] at (2.15,0) [anchor=west] {652};
\draw [fill=white] (4,0) circle (3.0pt) node[anchor=west] {~680};
\end{tikzpicture} 
\end{tabular}
\caption{IEEE 13 node radial distribution test feeder with parallel lines. Black lines denote existing lines and red dashed lines denote possible expansions.}
\label{P2-Fig:IEEE13}
\end{figure}

The original IEEE 13 node test feeder has 13 nodes and 12 lines (black solid lines in Figure \ref{P2-Fig:IEEE13}). For this paper, the network is modified  as follows. Continuous resources ($C_1$ through $C_5$) can be installed at 
nodes $611$ and $675$. Discrete resources ($D_1$ through $D_5$) can be installed at nodes $645$ and $650$.  Both types of resources can be installed at node $652$. The load for this system is based on data from a New Mexico distribution utility. Demand is added at all nodes and for all time-steps except for nodes $633$, $650$, $680$, $684$, and $692$. These nodes have zero demand during the entire design horizon. The installation and operational costs for all resources are in Table \ref{P2-Tab:TechChar}.

\begin{table}[ht!]
\caption{Characteristics of the technology options used in the IEEE 13 Network.}
\scriptsize
\begin{center}
\begin{tabular}{ p{0.6cm} p{1.0cm} p{1.0cm} p{2.0cm} p{1.8cm}}
\toprule
Tech Type & Fixed Cost & Variable Cost & Operational Cost  $\bs{a}P^2+\bs{b}P+\bs{c}$ & Rated Power (Max, Min)\\
&(\$)&(\$/KW)&(\$)& (KW)\\
\cmidrule{1-5}
$*C_1$ & 100,000 & 300 & $10P^2 + 5P + 2$ & ($100$, $-100$)\\
$C_2$ & 200,000 & 250 & $20P^2 + 10P + 4$ & ($100$,  $-100$) \\
$C_3$ & 250,000 & 200 & $30P^2 + 15P + 8$ & ($100$,  $-100$) \\
$C_4$ & 300,000 & 150 & $40P^2 + 20P + 10$ & ($100$,  $-100$) \\
$C_5$ & 350,000 & 100 & $50P^2 + 25P + 5$ & ($100$,  $-100$) \\
$D_1$ & 200,000 & 0 & $50P^2 + 25P + 6$ & ($250$, $-250$)\\
$D_2$ & 100,000 & 0 & $40P^2 + 20P + 5$ & ($275$, $-250$)\\
$D_3$ & 250,000 & 0 & $30P^2 + 15P + 4$ & ($300$, $-250$)\\
$D_4$ & 300,000 & 0 & $20P^2 + 10P + 3$ & ($225$, $-250$)\\
$D_5$ & 350,000 & 0 & $10P^2 + 5P + 2$ & ($200$, $-250$)\\
\cmidrule{1-5}
\multicolumn{5}{l}{$*$ indicates storage devices and the units are in \$/KVA}\\
\bottomrule
\end{tabular}\label{P2-Tab:TechChar}
\end{center}
\end{table}


Expansion decisions for this network include parallel lines for all 12 existing lines and new lines between nodes $611$ and $646$ (parallel and new lines are marked as dotted red lines in Figure \ref{P2-Fig:IEEE13}).  The cost of installing parallel lines and new lines is \$1000 per line. Physical characteristics of the lines are provided in \cite{chalil2017remote} and \cite{kersting2001radial}.  Here, $|\mathcal{T}| = 96$ and models 24 hours in 15 minute increments. This network has 18 possible generator contingencies and 25 possible line contingencies. 

\subsubsection*{Recommended solution for 96 design horizon problem for IEEE 13 network}
In this model, the optimal solution includes the installation of $D_2$ generators at nodes $650$ and $652$. The optimal solution also includes parallel lines between nodes $632$ -- $633$, $633$ -- $634$, $671$ -- $692$, $692$ -- $675$, and $611$ -- $646$.  The total installation and operational costs for 96 design horizon for this model is 
$\sim \$ 357k$
and is N-1 secure.

\subsubsection*{Solution time}
\begin{figure}[ht!]
\centering 
\subfigure[Solution times]{%
\label{P2-Fig:ResultsBA-SBD-RH-ST}%
\begin{tikzpicture}[scale=0.5, every node/.style={transform shape}]
\begin{semilogyaxis}[
ybar,
ymin = 10,
ymax = 70000,
xmin = 5,
xmax = 96,
enlarge y limits={0.25,upper},
enlarge x limits=0.10,
legend style={at={(0.45,0.95)},
legend style={draw=none},
legend style={/tikz/every even column/.append style={column sep=0.4cm}},
anchor=north,legend columns=4},
symbolic x coords={5,10,15,20,50,96},
xtick={5,10,15,20,50,96},
ylabel={Time (seconds)},
xlabel={Design Horizon},
point meta=rawy,
x tick label style={rotate=0,anchor=north},
    nodes near coords,
    every node near coord/.append style={font=\large, inner sep=2pt,anchor=west,rotate=90},  
    ticklabel style = {font=\large},    
]
\addplot coordinates {(5,336) (10,917) (15,2521) (20,4480) (50,12587) (96,52274)};
\addplot coordinates {(5,1022) (10,531) (15,1191) (20,1609) (50,7414) (96,12219)};
\addplot coordinates {(5,265) (10,4033) (15,4130) (20,5313) (50,2390) (96,3245)};
\addplot coordinates {(5,3193) (10,499) (15,517) (20,521) (50,667) (96,765)};
\legend{Base Algorithm, SBD, RH, SBD+RH}
\end{semilogyaxis}
\end{tikzpicture}
}%
\qquad
\subfigure[Objective value]{%
\label{P2-Fig:ResultsBA-SBD-RH-OF}%
\begin{tikzpicture}[scale=0.5, every node/.style={transform shape}]
\begin{axis}[
ymin = 200000,
ymax = 360000,
symbolic x coords={5,10,15,20,50,96},
xtick={5,10,15,20,50,96},
ytick={200000,250000,300000,350000,400000},
ylabel near ticks, yticklabel pos=left,
enlarge y limits={0.25,upper},
enlarge x limits=0.15,
legend style={at={(0.45,0.95)},
legend style={draw=none},
legend style={/tikz/every even column/.append style={column sep=0.4cm}},
anchor=north,legend columns=4},
ylabel={Objective value (\$)},
xlabel={Design Horizon},
point meta=rawy,
    x tick label style={rotate=0,anchor=north},
    nodes near coords,
    every node near coord/.append style={font=\large, inner sep=2pt,anchor=north west,rotate=0},    
]
\addplot coordinates {(5,217650) (10,224388) (15,229423) (20,234951) (50,268264)(96,357431)};
\addplot coordinates {(5,217650) (10,224388) (15,229423) (20,234951) (50,268264)(96,357431)};
\addplot coordinates {(5,217650) (10,224388) (15,229423) (20,234951) (50,268264)(96,357431)};
\addplot coordinates {(5,217650) (10,224388) (15,229423) (20,234951) (50,268264)(96,357431)};
\legend{Base Algorithm, SBD, RH, SBD+RH}
\end{axis}
\end{tikzpicture}
}%
\caption{Solution times and solution quality for the IEEE 13 case. RH refers to the rolling-horizon algorithm where the base algorithm is used to solve sub problems.  SBD+RH refers to the rolling-horizon algorithm where SBD is used to solve sub problems.}
\label{P2-Fig:ResultsBA-SBD-RH}
\end{figure}
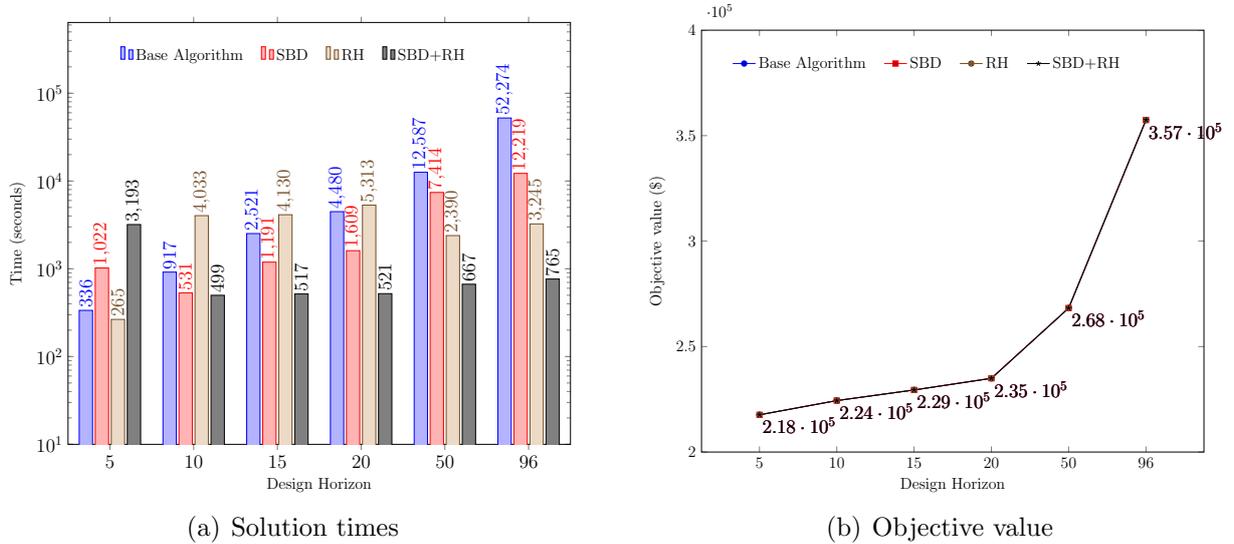

Figure \ref{P2-Fig:ResultsBA-SBD-RH} evaluates the efficiency and effectiveness of our RH approach. This figure shows the computation time and solution quality of the two exact methods and RH for design horizons of 5, 10, 15, 20, 50, and 96.  Each algorithm had a time limit of 24 hours. In all cases, the RH solution matches the solution found by the exact methods (Base algorithm, SBD). The SBD+RH was able to solve the 96 design horizon problem in 765 seconds ($\approx$ 68 times faster than the base algorithm). These results indicate that RH algorithm is computationally more efficient compared to the exact methods since it provides solutions of the same quality in a short time. However, it is important to stress that the RH is a heuristic \cite{silvente2015rolling}.

\subsubsection{Sensitivity analysis using IEEE 13 node test feeder}

\textit{Line installation costs}: In Table \ref{P2-Tab:SenseLineCosts}, we show results that indicate how the solution changes as the cost of building lines increases.
In this first case (1), the cost of building lines is \$1000. Here, the solution is to build 2 discrete generators and 5 lines. In the second, third, and fourth cases, the cost of building lines is increased to \$10,000, \$100,000, and \$1,000,000 respectively. On the last two cases, the solution is to build fewer lines and build a more expensive generator at node 675 (instead of node 645) to support all the contingencies.

\begin{table}[ht!]
\caption{Results that describe how the solution to the IEEE 13 node microgrid problem changes as the cost of adding lines increases.}
\scriptsize
\begin{center}
\begin{tabular}{ll|l|l}
\toprule
Case & Line Cost  & \multicolumn{2}{c}{Model Install Decisions} \\
\cmidrule{3-4}
 & (\$ per line)   & Generators & Lines (From -- To)  \\
\cmidrule{1-4}
1 & 1,000  & $D_2$ at Node 650 & (646 -- 611), (632 -- 633), (633 -- 634), \\
  &       & $D_2$ at Node 652 & (671 -- 692), (692 -- 675)\\
\cmidrule{1-4}
2 & 10,000  & $D_2$ at Node 650 & (646 -- 611), (632 -- 633), (633 -- 634), \\
  &       & $D_2$ at Node 652 & (671 -- 692), (692 -- 675)\\
\cmidrule{1-4}
3 & 100,000  & $D_2$ at Node 652 & (646 -- 611), (632 -- 633), (633 -- 634) \\
  &           & $C_2$ at Node 675 & \\
\cmidrule{1-4}
4 & 1,000,000  & $D_2$ at Node 652 & (646 -- 611), (632 -- 633), (633 -- 634) \\
  &           & $C_2$ at Node 675  & \\
\bottomrule
\end{tabular}\label{P2-Tab:SenseLineCosts}
\end{center}
\end{table}

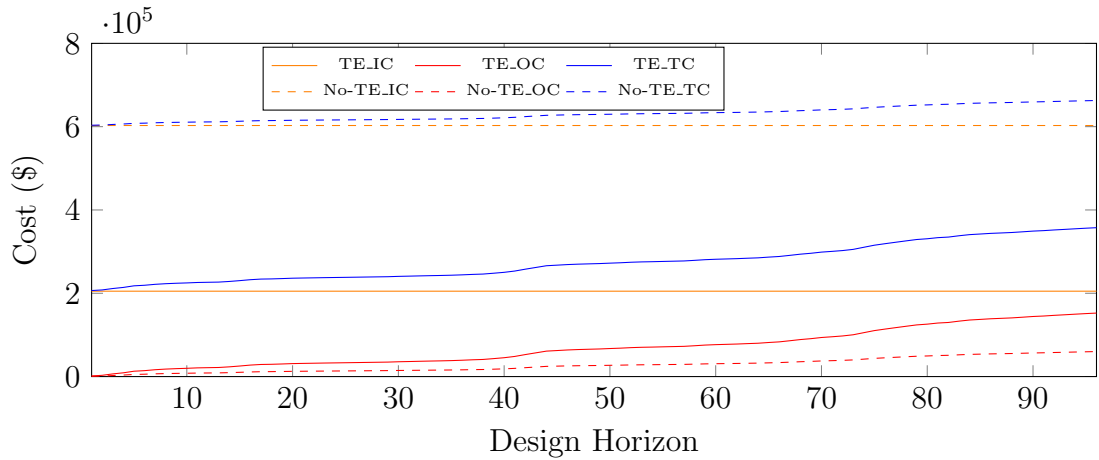
\begin{figure}[ht!]
\begin{tikzpicture}
\begin{axis}[
    enlargelimits=false,
    width=.9\textwidth,
    height=6cm,
    xlabel=Design Horizon,
	ylabel=Cost (\$),
	legend style={at={(0.4,0.99)},
		anchor=north,legend columns=3},
ymin = 0,
ymax = 800000,		
]
\addplot[color=orange] table[meta=Design_Horizon] {Data/12-ITGHC-IC96.tex};	
\addlegendentry{\tiny{TE\_IC}}
\addplot[color=red] table[meta=Design_Horizon] {Data/12-ITGHC-OC96.tex};	
\addlegendentry{\tiny{TE\_OC}}
\addplot[color=blue] table[meta=Design_Horizon] {Data/12-ITGHC-TC96.tex};	
\addlegendentry{\tiny{TE\_TC}}
\addplot[color=orange,dashed] table[meta=Design_Horizon] {Data/13-ITNLGHC-IC96.tex};	
\addlegendentry{\tiny{No-TE\_IC}}
\addplot[color=red,dashed] table[meta=Design_Horizon] {Data/13-ITNLGHC-OC96.tex};	
\addlegendentry{\tiny{No-TE\_OC}}
\addplot[color=blue,dashed] table[meta=Design_Horizon] {Data/13-ITNLGHC-TC96.tex};	
\addlegendentry{\tiny{No-TE\_TC}}
\end{axis}
\end{tikzpicture}
\caption{Operating cost and Installation cost IEEE 13 network.}
\label{P2-Fig:CumCostIGHCvsIGHCNL}
\end{figure}

\textit{Topology options}: We test the model with topology expansion options (TE) and without topology expansion options (No-TE). Table \ref{P2-Tab:SenseTopOpt} describes the optimal solution for IEEE 13 case with a time horizon of 96 for TE and No-TE. The solution for the No-TE model has load shedding at node $634$ if lines $632$ -- $633$ or $633$ -- $634$ fails. Figure \ref{P2-Fig:CumCostIGHCvsIGHCNL} shows how the costs of TE and No-TE change as the design horizon increases.
The costs are separated by installation costs (IC), operation cost (OC), and total costs (TC).

\begin{table}[ht!]
\caption{Sensitivity analysis on topology expansion options for IEEE 13 node test feeder}
\scriptsize
\begin{center}
\begin{tabular}{l|rl|rl}
\toprule
& \multicolumn{2}{c|}{With topology expansion} & \multicolumn{2}{c}{No topology expansion}\\
Cost type &  \multicolumn{1}{c}{Cost value}  & Optimal solution & \multicolumn{1}{c}{Cost value}  & Optimal solution \\
\cmidrule{1-5}
\multirow{2}{*}{Total generator installation cost} & \multirow{2}{*}{\$200,000.00}  & \multirow{2}{*}{$D_2$ at $650$ and $652$} & \multirow{2}{*}{\$ 602,735.45} & $C_2$ at $611$ and $675$\\
& & & & $D_2$ at $645$ and $652$\\
\cmidrule{1-5}
\multirow{3}{*}{Total line installation cost} & \multirow{3}{*}{\$5,000.00} &$646$ -- $611$ &\multirow{3}{*}{\$0.00} &\multirow{3}{*}{N/A}\\
& &$632$ -- $633$, $633$ -- $634$ & &\\
& &$671$ -- $692$, $692$ -- $675$ & &\\
\cmidrule{1-5}
Total operation cost (96 design horizon)& \$152,430.65 & &\$60,039.00 & \\
Total cost (96 design horizon) & \$357,430.65 & & \$662,774.45 & \\
& \multicolumn{2}{l|}{N-1 Secure on all nodes} & \multicolumn{2}{l}{No N-1 Security on node $634$}\\
\bottomrule
\end{tabular}\label{P2-Tab:SenseTopOpt}
\end{center}
\end{table}

\begin{figure}[ht!]
\centering
\begin{tabular}{cc}
\begin{tikzpicture}[scale=1]
\draw (0,2) -- (8,2);
\draw (0,4) -- (8,4);
\draw (2,0) -- (2,2);
\draw (4,0) -- (4,6);
\draw[red,dashed] (0,1.9) -- (8,1.9);
\draw[red,dashed] (0,3.9) -- (8,3.9);
\draw[red,dashed] (1.9,0) -- (1.9,1.9);
\draw[red,dashed] (3.9,0.1) -- (3.9,6);
\draw[red,dashed] (0,2) -- (0,3.9);
\draw [fill=white] (0,4) circle (3.0pt) node[anchor=south] {646} ;
\draw [fill=blue!50,draw=blue] (-0.15,4.1)--(0.15,4.1)--(0,3.8) -- cycle;
\filldraw[fill=magenta!50,draw=magenta] (2,4) ellipse (1.5pt and 3pt) node[anchor=south] {645};
\filldraw[fill=magenta!50,draw=magenta] (1.8,4) ellipse (1.5pt and 3pt);
\filldraw[fill=magenta!50,draw=magenta] (2.2,4) ellipse (1.5pt and 3pt);
\draw [fill=white] (4,4) circle (3.0pt) node[anchor=south west] {632};
\draw [fill=blue!50,draw=blue] (3.85,4.1)--(4.15,4.1)--(4,3.8) -- cycle;
\draw [fill=white] (6,4) circle (3.0pt) node[anchor=south] {633};
\draw [fill=white] (8,4) circle (3.0pt) node[anchor=south] {634};
\draw [fill=blue!50,draw=blue] (7.85,4.1)--(8.15,4.1)--(8,3.8) -- cycle;
\filldraw[fill=magenta!50,draw=magenta] (4,6) ellipse (3pt and 1.5pt) node[anchor=south] {650};
\filldraw[fill=red!50, draw=red] (-0.1,2.1) rectangle (0.1,1.9) node[anchor=north] at (0,1.9) {611};
\draw [fill=white] (2,2) circle (3.0pt) node[anchor=south] {684};
\draw [fill=white] (4,2) circle (3.0pt) node[anchor=south west] {671}; 
\draw [fill=blue!50,draw=blue] (3.85,2.1)--(4.15,2.1)--(4,1.8) -- cycle;
\draw [fill=white] (6,2) circle (3.0pt) node[anchor=south] {692};
\filldraw[fill=red!50, draw=red]  (7.9,1.9) rectangle (8.1,2.1)  node[anchor=south] {675};
\draw [fill=blue!50,draw=blue] (1.85,0)--(2.15,0)--(2,-0.25)--cycle;
\node [blue, fill=white] at (2.15,0) [anchor=west] {652};
\draw [fill=white] (4,0) circle (3.0pt) node[anchor=west] {~680};
\end{tikzpicture} 
\end{tabular}
\caption{Modified IEEE 13 Node Network from  Figure \ref{P2-Fig:IEEE13}}
\label{P2-Fig:IEEE13-C2}
\end{figure}
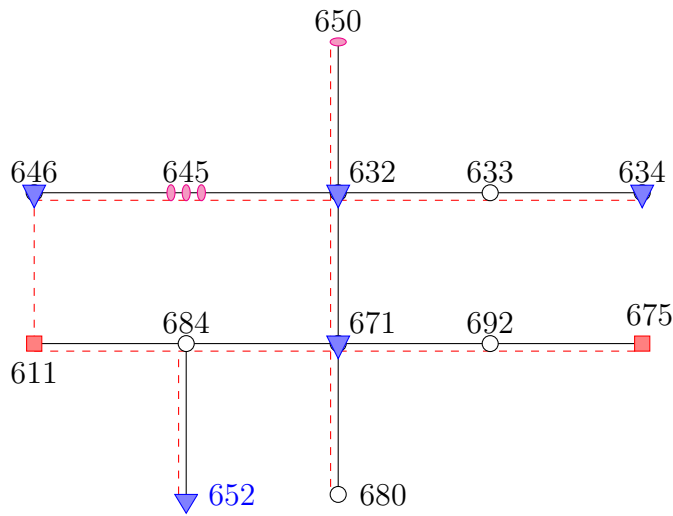

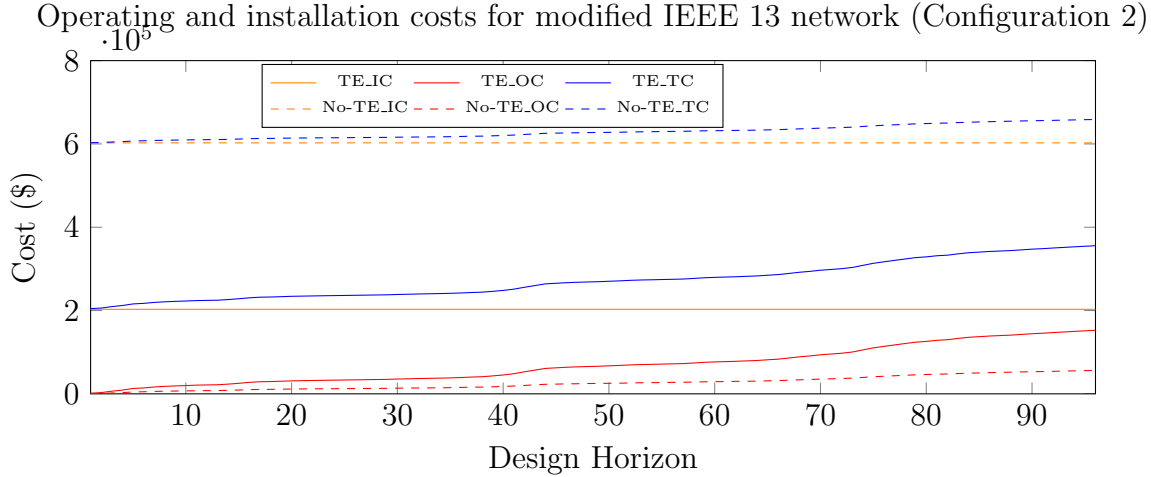
\begin{figure}[ht!]
\begin{tikzpicture}
\begin{axis}[
    title=Operating and installation costs for modified IEEE 13 network (Configuration 2),
    enlargelimits=false,
    width=.9\textwidth,
    height=6cm,
    xlabel=Design Horizon,
	ylabel=Cost (\$),
	legend style={at={(0.4,0.99)},
		anchor=north,legend columns=3},
ymin = 0,
ymax = 800000,		
]
\addplot[color=orange] table[meta=Design_Horizon] {Data/14-ITAGHC-IC96.tex};	
\addlegendentry{\tiny{TE\_IC}}
\addplot[color=red] table[meta=Design_Horizon] {Data/14-ITAGHC-OC96.tex};	
\addlegendentry{\tiny{TE\_OC}}
\addplot[color=blue] table[meta=Design_Horizon] {Data/14-ITAGHC-TC96.tex};	
\addlegendentry{\tiny{TE\_TC}}
\addplot[color=orange,dashed] table[meta=Design_Horizon] {Data/15-ITAGNLHC-IC96.tex};	
\addlegendentry{\tiny{No-TE\_IC}}
\addplot[color=red,dashed] table[meta=Design_Horizon] {Data/15-ITAGNLHC-OC96.tex};	
\addlegendentry{\tiny{No-TE\_OC}}
\addplot[color=blue,dashed] table[meta=Design_Horizon] {Data/15-ITAGNLHC-TC96.tex};	
\addlegendentry{\tiny{No-TE\_TC}}
\end{axis}
\end{tikzpicture}
\caption{Operating cost and Installation cost for IEEE 13 Network with configuration 2 in Figure \ref{P2-Fig:IEEE13-C2}.}
\label{P2-Fig:CumCostIGAGHCvsIGAGHCNL}
\end{figure}

We also modify the network in Figure \ref{P2-Fig:IEEE13} as shown in Figure \ref{P2-Fig:IEEE13-C2} with the option to install generators at all demand nodes to compare the impacts of including topology expansion options. this also guarantees power supply during N-1 contingencies at all demand nodes. In this model, with TE options, the optimal solution includes the installation of $D_2$ generators at nodes $634$ and $652$ and parallel lines between nodes $671$ -- $692$, $692$ -- $675$ and $611$ -- $646$. The total installation and operation costs for this modified network is 
$\sim \$ 355k$.
In contrast, the optimal solution for No-TE includes installation of $D_2$ generators at nodes $634$, $646$, and $652$ and $C_2$ generators at nodes $611$ and $675$. The total costs for No-TE is 
$\sim \$ 658K$.
Figure \ref{P2-Fig:CumCostIGAGHCvsIGAGHCNL} provides cumulative costs with the progression of the design horizon for both TE and No-TE options for configuration 2. 

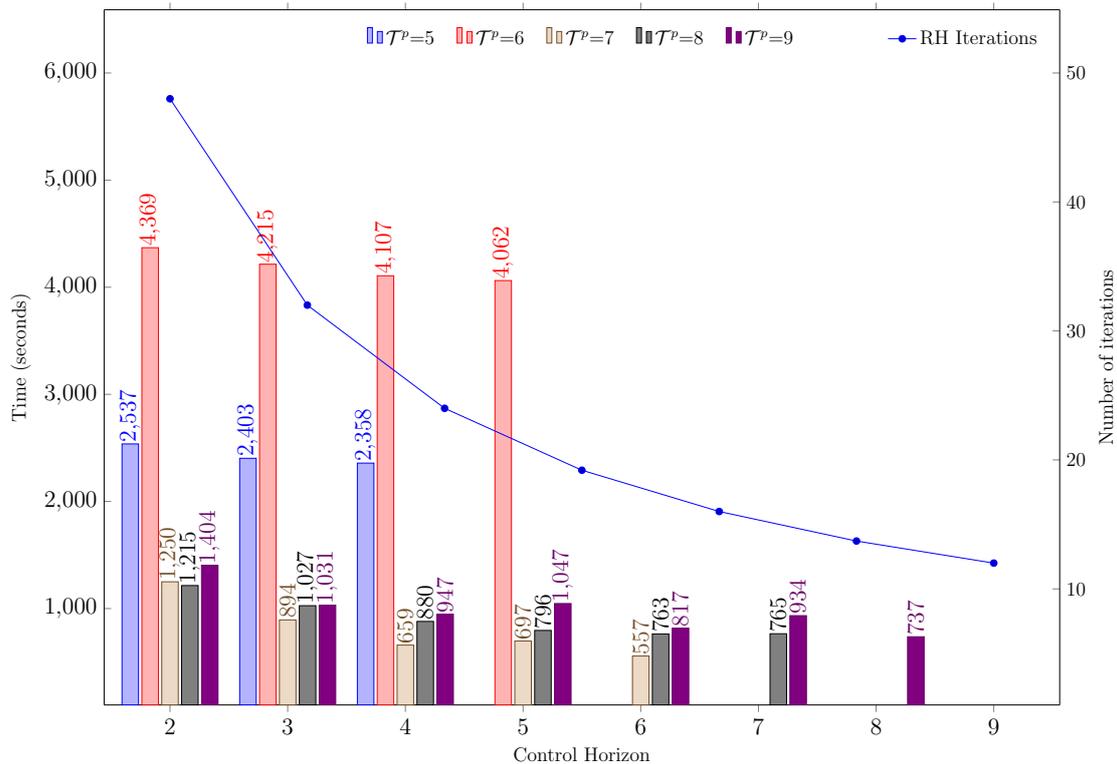
\begin{figure}[ht!]
\centering 
\begin{tikzpicture}[scale=0.62, every node/.style={transform shape}]
\begin{axis}[
ybar,
axis y line*=left,
 ylabel=y-axis 1,
ymin = 100,
ymax = 6000,
xmin = 2,
xmax = 9,
width=\textheight,
height=\textwidth,
enlarge y limits={0.10,upper},
enlarge x limits=0.08,
legend style={at={(0.5,0.98)},
legend style={draw=none},
legend style={/tikz/every even column/.append style={column sep=0.4cm}},
anchor=north,legend columns=-1},
symbolic x coords={2,3,4,5,6,7,8,9},
xtick={2,3,4,5,6,7,8,9},
ytick={0,1000,2000,3000,4000,5000,6000},
ylabel={Time (seconds)},
xlabel={Control Horizon},
point meta=rawy,
x tick label style={rotate=0,anchor=north},
    nodes near coords,
    every node near coord/.append style={font=\large, inner sep=1pt,anchor=west,rotate=90},  
    ticklabel style = {font=\large}, 
]
\addplot coordinates {(2,2537) (3,2403) (4,2358)};
\addplot coordinates {(2,4369) (3,4215) (4,4107) (5,4062)};
\addplot coordinates {(2,1250) (3,894) (4,659) (5,697) (6,557)};
\addplot coordinates {(2,1215) (3,1027) (4,880) (5,796) (6,763) (7,765)};
\addplot coordinates {(2,1404) (3,1031) (4,947) (5,1047) (6,817) (7,934) (8,737)};
\legend{$\mathcal{T}^p$=5, $\mathcal{T}^p$=6, $\mathcal{T}^p$=7, $\mathcal{T}^p$=8,$\mathcal{T}^p$=9}
\end{axis}
\begin{axis}[
axis y line*=right,
 ylabel=y-axis 2,
ymin = 1,
ymax = 50,
xmin = 2,
xmax = 8,
width=\textheight,
height=\textwidth,
enlarge y limits={0.10,upper},
enlarge x limits=0.08,
legend style={at={(0.9,0.98)},
legend style={draw=none},
anchor=north,legend columns=-1},
symbolic x coords={2,3,4,5,6,7,8},
ytick={10,20,30,40,50},
ylabel={Number of iterations},
point meta=rawy,
axis x line=none
]
\addplot coordinates {(2,48) (3,32) (4,24) (5,19.2) (6,16) (7,13.71) (8,12)};
\legend{RH Iterations}
\end{axis}
\end{tikzpicture}
\caption{Parameter experiment using various $\mathcal{T}^p$ and $\mathcal{T}^c$}%
\label{P2-Fig:ResultsRHCHWSabc}
\end{figure}

\textit{Rolling-horizon parameters}: The choices of $\mathcal{T}^p$ and $\mathcal{T}^c$ can have an impact on solution quality.  Results that vary these parameters are presented in Figure \ref{P2-Fig:ResultsRHCHWSabc}. In this case the solution remains the same for all experiments. However, the solution time varies significantly. The structure of this distribution system forces solutions to add significant redundancy to satisfy N-1 security constraints. Thus, the RH with small $\mathcal{T}^p$ and $\mathcal{T}^c$ has sufficient information to make decisions that are of high quality for the entire planning horizon.
The time variation is due to two characteristics of the RH algorithm: number of sub-problems and length of $\mathcal{T}^p$. The number of sub-problems is a function of the ratio $\frac{|\mathcal{T}|}{|\mathcal{T}^c|}$. For example, a 96 design horizon ($|\mathcal{T}|$) problem with control horizon ($|\mathcal{T}^c|$) of 4 requires 24 iterations. As the number of iterations is reduced, it is faster to solve the model. But at the same time, as the length of prediction horizon increases, it takes longer time to solve the initial iteration (Figure \ref{P2-Fig:ResultsBA-SBD-RH}). 

\textit{Demand changes}: We conduct a sensitivity analysis with respect to demand. We rerun the model by multiplying each of the base demands ($\bs{dp}_i^t$ and $\bs{dq}_i^t$) with the values 0.20, 0.40, 0.80, 1.20, and 1.40 to check if the model decisions changed with the change in demand. From the results, all the solutions recommended installation of two generators and five lines. The power generated by these generators changed to match the changing demands.


\subsection{Case study 2: Alaskan microgrid}

\begin{figure}[ht!]
\centering 
\subfigure[Configuration 1]{%
\label{P2-Fig:Nome}%
\includegraphics[width=0.46\textwidth]{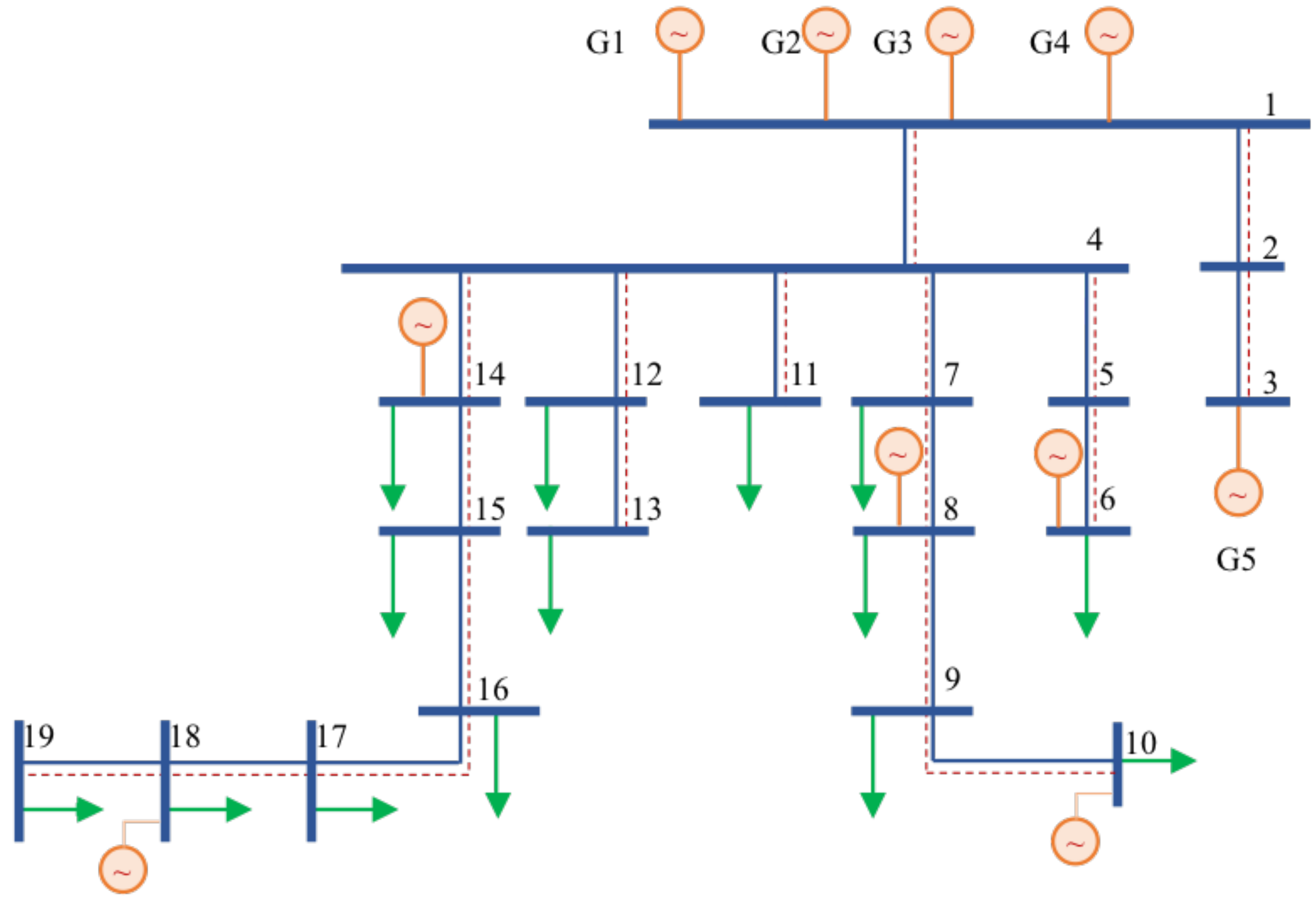}}%
\qquad
\subfigure[Configuration 2]{%
\label{P2-Fig:Nome1}%
\includegraphics[width=0.46\textwidth]{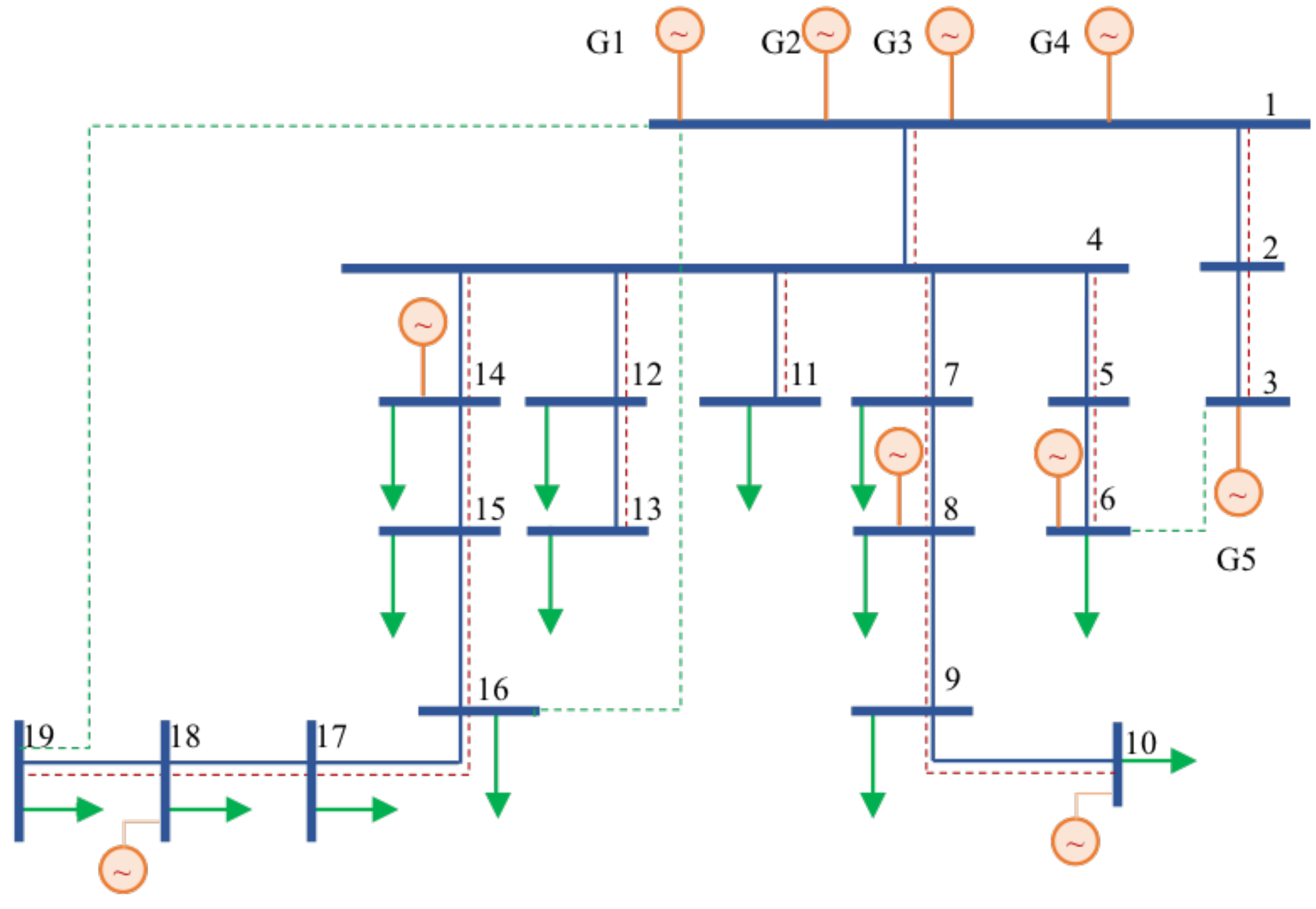}}%
\caption{Schematic diagram of a remote community in Alaska.}
\end{figure}

In this section, we test the performance of the RH algorithm on two variations of an Alaskan micogrid that has 19 nodes and 18 lines (Figure \ref{P2-Fig:Nome}).
In the first variation of this model (referred to as configuration 1), a single installation of a discrete resource ($D_1$ through $D_5$) is allowed at each node $6, 8, 10, 14,$ and $18$ (Table \ref{P2-Tab:TechCharNome}).
Demands for this system are based on data provided by the Alaskan distribution utility \cite{Mashayekh2016}.  Installation and operational costs are provided in Table \ref{P2-Tab:TechCharNome}. 
Table \ref{P2-Tab:LineCharNome} describes the specifications of the lines in this system.  This network has four existing generators at node $1$ and one existing generator at node $3$. All five existing generators are of type $D_1$. The capacity of the existing generators is modified so that the model is forced to build new generators. Parallel lines may be built anywhere in the system, provided  a line currently exists, at a cost of  \$1,000 per line. There are seven generator contingencies and 36 possible line contingencies in the network. The seven generator contingencies are due to five new generators, one existing generator node $3$, and one of the 5 generators from node $1$.  Altogether, there are 43 contingencies for this network. 

In the second version of the model (referred to as configuration 2), we modify the system as shown in Figure \ref{P2-Fig:Nome1}. In this configuration, options to install 3 new lines between nodes 1--19, 1--16, and 3--6 are added.

\begin{table}[ht!]
\caption{Characteristics of technology options for Case Study 2.}
\scriptsize
\begin{center}
\begin{tabular}{ p{1.2cm} p{0.8cm} p{2.2cm} p{1.9cm} }
\toprule
Tech Type & Fixed Cost & Operational Cost  $\bs{a}P^2+\bs{b}P+\bs{c}$ & Rated Power (Max, Min)\\
&(\$)&(\$)& (KW)\\
\cmidrule{1-4}
$D_1$ & $200,000$ & $50P^2 + 25P + 6$ & ($200$, $-200$)\\
$D_2$ & $500,000$ & $60P^2 + 20P + 5$ & ($1500$, $-1500$)\\
\bottomrule
\end{tabular}\label{P2-Tab:TechCharNome}
\end{center}
\end{table}

\begin{table}[ht!]
\caption{Line configuration for Case study 2.}
\scriptsize
\begin{center}
\begin{tabular}{lllll}
\toprule
ID & Resistance  & Reactance   & Thermal Limit & Lines \\
 &  pu & pu   & $MVA$ & (From node - To node) \\
\cmidrule{1-5}
A & 0  & 0.05 & 10000  & (1--2), (1--4), (5--6) \\
B  & 0.392921923 & 0.923131194 & 3422.532396   & (2--3)   \\
C  & 0.157168769 & 0.369252478 & 3422.532396   & (4--5)   \\
D  & 0.002854927 & 0.005210712 & 3782.798964   & (4--7), (7--8)  \\
E  & 0.019646096 & 0.04615656  & 3422.532396   & (8--9)   \\
F  & 0.039292192 & 0.092313119 & 3422.532396   & (9--10)  \\
G  & 0.314337539 & 0.738504956 & 3422.532396   & (4--11), (4--12)\\
H  & 0.1021597   & 0.240014111 & 3422.532396   & (12--13) \\
I  & 0.248685034 & 0.20405299  & 1585.172899   & (4--14)  \\
J  & 0.373027551 & 0.306079486 & 1585.172899   & (14--15) \\
K  & 0.062171258 & 0.051013248 & 1585.172899   & (15--16) \\
L  & 0.497370068 & 0.408105981 & 1585.172899   & (16--17), (17--18)\\
M  & 1.243425169 & 1.020264952 & 1585.172899   & (18--19)  \\
N$^*$  & 1.0 & 1.0 & 3782.798964   & (1--19), (1--16), (3--6)  \\
\cmidrule{1-5}
\multicolumn{5}{l}{$^*$ indicates new line characteristics for configuration 2 from Figure \ref{P2-Fig:Nome1}}  \\
\bottomrule
\end{tabular}\label{P2-Tab:LineCharNome}
\end{center}
\end{table}

\subsubsection*{Recommended solution for 96 design horizon for Alaskan microgrid network}
\textit{Configuration 1:} The solution for this microgrid installs generators of type $D2$ at nodes $6$, $8$, $10$, $14$ and $18$. The solution also installs parallel lines in all locations (18) to support N-1 security. The total installation cost is $\sim \$ 2,518K$ 
and the total operating cost is 
\$7,707,405,649 
$\sim \$ 7,707,405K$.
(Total cost = 
$\sim \$ 7,709,923K$
).

\textit{Configuration 2:} In configuration 2, the
solution to $\cal M$
also installs generators of type $D2$ at nodes $6$, $8$, $10$, $14$ and $18$. The solution installs 19 parallel lines except the lines between nodes (15--16) and (1--16). In comparison with configuration 1, 17 parallel lines and 2 new lines are installed to satisfy N-1 security constraints. However, the total operating costs is lower than the configuration 1. The total installation cost is \$2,519,000 and total operating cost is \$7,707,278,147 (Total cost = \$7,709,797,147).

\subsubsection*{Solution time}
\begin{figure}[ht!]
\centering 
\label{P2-Fig:ResultsBA-SBD-RH-Nome}
\subfigure[Solution times]{%
\label{P2-Fig:ResultsBA-SBD-RH-Nome-ST}%
\begin{tikzpicture}[scale=0.5, every node/.style={transform shape}]
\begin{semilogyaxis}[
ybar,
ymin = 10,
ymax = 70000,
xmin = 5,
xmax = 96,
enlarge y limits={0.25,upper},
enlarge x limits=0.10,
legend style={at={(0.35,0.95)},
legend style={draw=none},
legend style={/tikz/every even column/.append style={column sep=0.2cm}},
anchor=north,legend columns=4},
symbolic x coords={5,10,15,20,50,96},
xtick={5,10,15,20,50,96},
ylabel={Time (seconds)},
xlabel={Design Horizon},
point meta=rawy,
x tick label style={rotate=0,anchor=north},
    nodes near coords,
    every node near coord/.append style={font=\large, inner sep=2pt,anchor=west,rotate=90},  
    ticklabel style = {font=\large},    
]
\addplot coordinates {(5,819) (10,1767) (15,6380) (20,13513) (50,86400) (96,86400)};
\addplot coordinates {(5,859) (10,2399) (15,5478) (20,11141) (50,86400) (96,86400)};
\addplot coordinates {(5,308) (10,2595) (15,2601) (20,2622) (50,1928) (96,2048)};
\addplot coordinates {(5,624) (10,1746) (15,1756) (20,1800) (50,646) (96,877)};
\legend{Base Algorithm, SBD, RH, SBD+RH}
\end{semilogyaxis}
\end{tikzpicture}
}%
\qquad
\subfigure[Objective values]{%
\label{P2-Fig:ResultsBA-SBD-RH-Nome-OF}%
\begin{tikzpicture}[scale=0.5, every node/.style={transform shape}]
\begin{semilogyaxis}[
ymin = 3E8,
ymax = 8E9,
symbolic x coords={5,10,15,20,50,96},
xtick={5,10,15,20,50,96},
ylabel near ticks, yticklabel pos=left,
enlarge y limits={0.25,upper},
enlarge x limits=0.15,
legend style={at={(0.45,0.95)},
legend style={draw=none},
legend style={/tikz/every even column/.append style={column sep=0.4cm}},
anchor=north,legend columns=4},
ylabel={Objective value (\$)},
xlabel={Design Horizon},
point meta=rawy,
    x tick label style={rotate=0,anchor=north},
    nodes near coords,
    every node near coord/.append style={font=\large, inner sep=2pt,anchor=north west,rotate=0},    
]
\addplot coordinates {(5,4.29E8) (10,8.14E8) (15,1.17E9) (20,1.52E9) (50,1.02E10) };
\addplot coordinates {(5,4.29E8) (10,8.14E8) (15,1.17E9) (20,1.52E9) };
\addplot coordinates {(5,4.29E8) (10,8.14E8) (15,1.17E9) (20,1.52E9) (50,4.04E9)(96,7.71E9)};
\addplot coordinates {(5,4.29E8) (10,8.14E8) (15,1.17E9) (20,1.52E9) (50,4.04E9)(96,7.71E9)};
\legend{Base Algorithm, SBD, RH, SBD+RH}
\end{semilogyaxis}
\node[blue] at (10.5,8.7) {\normalsize{(81.2\% gap)}};
\end{tikzpicture}
}%
\caption{Solution times and solution quality for the Alaskan microgrid case. RH refers to the rolling-horizon algorithm where the base algorithm is used to solve sub problems.  SBD+RH refers to the rolling-horizon algorithm where SBD is used to solve sub problems. The base algorithm and SBD timed-out for 50 and 96 design horizon problems. The best feasible solution found for 50 design horizon was at 81.2\% optimality gap and for 96 design horizon, base algorithm and SBD failed to find any feasible solution even after 24 hours.}
\end{figure}
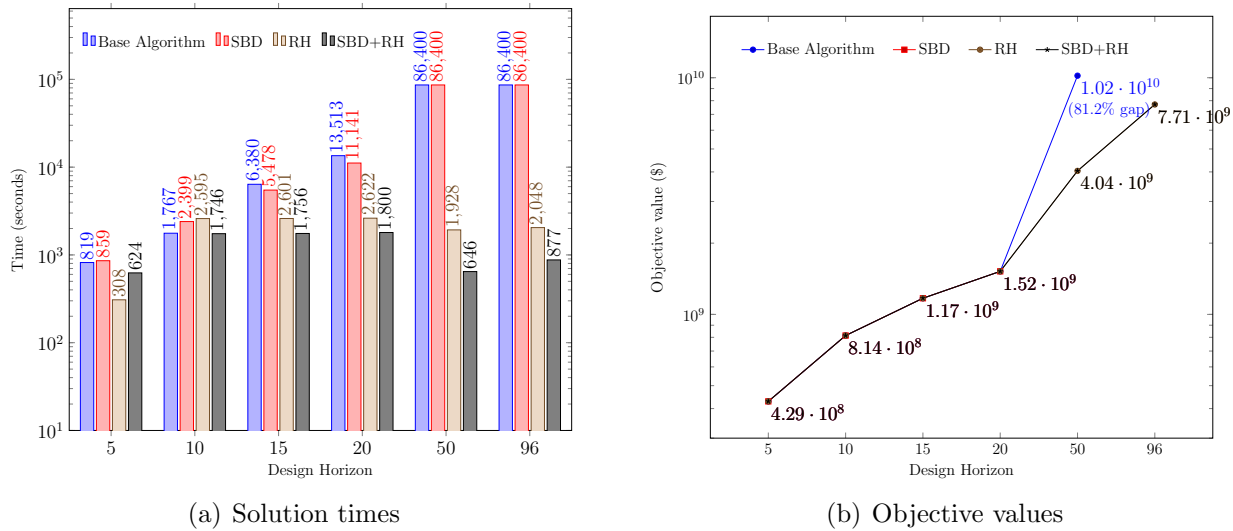

Figure \ref{P2-Fig:ResultsBA-SBD-RH-Nome} compares the solution time  for design horizons of 5, 10, 15, 20, 50, and 96. Here, all algorithms were terminated after 24 hours.
The two exact methods found the optimal solution for the case where $|\mathcal{T}| = 5, 10, 15, \textrm{ and } 20$.
The rolling-horizon algorithm with SBD is able to solve the 96 design horizon problem in 877 seconds.  On this problem, most of the contingencies must be added to $\mathcal{M}$, a situation that limits the effectiveness of SBD.

\section{Conclusions and Future Results}
In this paper, we develop a mathematical formulation for designing and operating remote off-grid microgrids with N-1 security constraints on generators and lines. We also present a rolling-horizon algorithm that efficiently solves these problems. The algorithm was tested on the IEEE 13 node system and a real Alaskan microgrid network. 
The focus of this paper was on extending existing microgrid design and operations models to handle topology expansion and line contingencies. Overall, the results suggest that topology expansions are needed to support low cost N-1 security in distribution systems. There remain a number of interesting future directions for this research. First, this model assumes that all generations and demands are deterministic. Future work should consider incorporating stochastic renewable resources, such as, wind and solar.  Introduction of wind and solar can result in various issues like line overheating and insufficient generation capacity and can be solved using probabilistic chance constrained approaches as used in \cite{sundar2016unit,bienstock2014chance}. Second, other solution algorithms should be explored to improve the scalabilty of solving ACIRMP, in both the size of the network and the length of the time horizons. Further, since rolling-horizon-based algorithms yield sub-optimal solutions, one can consider applying hierarchical rolling-horizon approaches to revisit the decisions made at earlier time windows \cite{Zavala2016,sai2018pscc}.

\paragraph{\textbf{Acknowledgements}} This work is partially funded by the Office of Electricity Delivery and Energy Reliability, Distributed Energy Program of the U.S. Department of Energy under Contract No. DE-AC02-05CH11231 and Work Order M615000466. We also want to thank Dan Ton from the Office of Electricity Delivery and Energy Reliability for supporting our work.
 Clemson University is acknowledged for generous allotment of compute time on Palmetto cluster. The Center for Nonlinear Studies at LANL also supported this work.  



\bibliographystyle{elsarticle-num}

\bibliography{bib_paper2.bib}

\end{document}